\documentclass [12pt]{article}
\usepackage{amssymb,amsmath, amsthm, amscd,ifthen}
\usepackage[T2A]{fontenc}
\usepackage[cp1251]{inputenc}
\usepackage[russian,english]{babel}
\usepackage[dvips]{graphicx}
\usepackage{indentfirst}
\newtheorem{Th}{Theorem}
\newtheorem{Lem}{Lemma}
\newtheorem{Stat}{Statement}

\newtheorem{cj}{Conjecture}

\newtheorem{cor}{Corollary}

\def\dist{\mathrm {dist}}

\title{The distribution of second degrees in the Buckley--Osthus random graph model}

\author{Andrey Kupavskii\footnote{Yandex, Department of theoretical and practical research; Department of Discrete Mathematics, Moscow Institute of Physics and Technology; Moscow State University, Faculty of Mechanics and Mathematics,
Department of Number Theory.},
Liudmila Ostroumova\footnote{Yandex, Department of theoretical and practical research; Moscow State University, Mechanics and Mathematics Faculty, Department of Mathematical Statistics and Random Processes.},
Dmitriy Shabanov\footnote{Yandex, Department of theoretical and practical research; Department of Discrete Mathematics, Moscow Institute of Physics and Technology; Moscow State University, Faculty of Mechanics and Mathematics,
Department of Probability Theory.}, Prasad Tetali\footnote{School of Mathematics and School of Computer Science, Georgia Institute of Technology, Atlanta, GA, USA}}

\date{}

\renewcommand{\le}{\leqslant}
\renewcommand{\ge}{\geqslant}

\newcommand{\B}{{\rm B}}

\begin{document}

\newcommand{\MExpect}{\mathsf{E}}
\newenvironment{Proof}{\noindent{\it Proof.\,}}{\hfill$\Box$}

\maketitle

\begin{abstract}
In this paper we consider a well-known generalization of the Barab\'asi and Albert preferential attachment model --  the Buckley-Osthus model. Buckley and Osthus proved that in this model the degree sequence has a power law distribution.
As a natural (and arguably more interesting) next step, we study the  second degrees of vertices.
Roughly speaking, the second degree of a vertex is the number of vertices at  distance two from this vertex. The distribution of second degrees is of interest because it
is a good approximation of PageRank, where the importance of a vertex is measured by taking  into account the popularity of its neighbors.

We prove that the second degrees also obey a power law. More precisely, we estimate the expectation of the number of vertices with the second degree greater than or equal to $k$ and prove the concentration of this random variable around its expectation using the now-famous Talagrand's concentration inequality over product spaces.
As far as we know this is the only application of Talagrand's inequality to random web graphs, where the  (preferential attachment) edges are not defined over a product distribution, making the application nontrivial, and requiring certain novelty.
\end{abstract}
\vskip+1cm

\section{Introduction}

In this paper we consider some properties of random graphs.
The standard random graph model $G(n,m)$ was introduced by Erd\H{o}s and R\'enyi in \cite{ER}.
In this model we randomly choose one graph from all graphs with $n$ vertices and $m$ edges.
The similar model $G(n,p)$ was suggested by Gilbert in \cite{G}:
we have an $n$-vertex set and we join each pair of vertices independently with probability $0<p<1$.
Many papers deal with the classical models.
Fundamental results about these models can be found in \cite{B}, \cite{K}, \cite{JLR}.

Recently there has been an increasing interest in modeling complex real-world networks.
It is well understood that real structures differ from standard random graphs. Many models of real-world networks and main results can be found in \cite{Res}.
For example, a basic characteristic of random graphs is their degree sequence.
In many real-world structures the degree sequence obeys a power law distribution.
However, standard random graph models do not have this property.

In 1999, Barab\'asi and Albert \cite{BA} suggested the so-called preferential attachment model that has a desired degree distribution.
Later Bollob\'as and Riordan \cite{BR2} gave a more precise definition of this model.
In this model the probability that a new vertex is connected to some previous vertex $v$
is proportional to the degree of $v$.
Bollob\'as and Riordan also proved that the degree sequence has a power law distribution
with exponent equal to $-3$.

Naturally one would not expect that this constant will suit all (or even most) of the real networks.
In order to make the model more flexible, two groups of authors (see \cite{Dor} and \cite{Drin}) proposed to add one more parameter --- an
``initial attractiveness'' of a node which is a positive constant that does not
depend on the degree.
In \cite{BO},  Buckley and Osthus gave an explicit construction of that model.

Many papers deal with different variations of preferential attachment.
We mention here the paper by Rudas, T\'oth and Valko (see \cite{Trees}).
The authors consider quite a generic model of a random tree and prove some
interesting results concerning a neighborhood structure of a random vertex.
Also one can find a neighborhood analysis in preferential attachment models
in the preprint \cite{Limits}
on the weak graph limit.

This paper deals with the Buckley--Osthus model, which we now describe.
Let $n$ be a number of vertices in our graph,
$m\in \mathbb{N}$ and $a\in\mathbb{R}_+$ be fixed parameters.

We begin with the case $m=1$.
We inductively construct a random graph $H_{a,1}^n$.
The graph  $H_{a,1}^1$ consists of one vertex and one loop (we can also start with $H_{a,1}^0$, which is the empty graph).
Assume that we have already constructed the graph $H_{a,1}^{t-1}$.
At the next step we add one vertex $t$ and one edge between vertices $t$ and $i$,
where $i$ is chosen randomly with
$$
\Prob(i=s) = \
\begin{cases}
\frac{d_{H_{a,1}^{t-1}}(s) - 1 + a}{(a+1)t-1} & \text{if } 1 \le s \le t-1,\cr \noalign{\smallskip}
\frac{a}{(a+1)t-1} &\text{if }  s=t. \cr
\end{cases}
$$
Here $d_{H_{a,1}^t}(s)$ is the degree of the vertex $s$ in $H_{a,1}^t$.
We will also use the notation $d(s):=d_{H_{a,1}^n}(s)$.

To construct $H_{a,m}^n$ with $m>1$ we start from  $H_{a,1}^{mn}$.
Then we identify the vertices $1, \dots, m$ to form the first vertex; we
identify the vertices $m+1, \dots, 2m$ to form the second vertex, etc.
As for the edges, if  the edge $e$ connects vertices $im+k$ and $jm+l, 1\le k,l\le m,$ in the graph $H_{a,1}^{mn}$ then we draw an edge $e'$ between vertices $i+1$ and $j+1$ in $H_{a,m}^n$. Note that we have a one-to-one correspondence between the edges of $H_{a,1}^{mn}$ and $H_{a,m}^n$, so there may be multiple edges (and multiple loops) between vertices in $H_{a,m}^n$.
Denote by $\mathfrak{H}_{a,m}^n$  the probability space of constructed graphs.

In \cite{BO} Buckley and Osthus proved that the degree sequence of $H_{a,m}^n$ has a power law
with exponent $-2-a$ if $a$ is a natural number.
Recently Grechnikov substantially improved this result.

\begin{Th}{\rm (Grechnikov, \cite{Gr})}
Let $a\in\mathbb{R}_+$. If $d = d(n) \ge m$ and $\psi(n) \to \infty$ as $n \to \infty$, then
$$
\left|R(d,n) - \frac{\B(d-m+ma,a+2)n} {\B(ma,a+1)}\right| \le \left(\sqrt{d^{-2-a}n} + d^{-1}\right)\psi(n)\,,
$$
with probability tending to $1$ as $n \to \infty$.
Here $R(d,n)$ is the number of vertices in $H_{a,m}^n$ with degree equal to $d$ and $\B(x,y)$ is the beta function.
\end{Th}


In this paper we consider the so-called second degrees of vertices. Roughly speaking, the second degree of a vertex is the number of vertices at distance two from the vertex. We prove that the number of vertices $Y_n(k)$ with the second degree at least $k$  decreases as $k^{-a},$ where $a$ is the initial attractiveness. This means that the distribution of second degrees obeys power law. To prove this we calculate the expectation of $Y_n(k)$ and show the concentration of this random variable around its expectation using Talagrand's inequality. The  application of this inequality is nontrivial, in particular, we have to redefine the probability space of the Buckley-Osthus graph so that we obtain a product probability space with a product measure. 
After that we use Talagrand's inequality in its general asymmetric form, which is essential.
Verifying the hypothesis  of Talagrand's theorem for the present purpose  turns out to be delicate, requiring us to introduce additional combinatorial constructions.

This paper is organized as follows.
In Section \ref{sec2} we give the main definitions and formulate the results.
In Sections \ref{sec3} and \ref{sec4} we prove the theorems stated in Section \ref{sec2}.

\section{Definitions and results}\label{sec2}
\subsection{Definitions}
In this paper we study the random graph $H_{a,1}^n$.
We shall write $ij \in H_{a,1}^n$ if  $H_{a,1}^n$ contains the edge $ij$;
we shall write $t \in H_{a,1}^n$ if $t$ is a vertex of $H_{a,1}^n$.
Given a vertex $ t \in H_{a,1}^n $,  the {\it second degree} of the vertex $t$ is
$$
d_2(t) = \#\{ij: i \not= t, j \not= t, it \in H_{a,1}^n, tj \in H_{a,1}^n\}.
$$
In other words, the second degree of $t$ is the number of edges adjacent to the neighbors of $t$
except for the edges adjacent to the vertex $t$.
We say that a vertex $t$ is a {\it $k$-vertex} if $d_2(t) \ge k$.

Let $M_n^1(d)$ be the expectation of the number of vertices with degree $d$ in $H_{a,1}^n$:
$$
M_n^1(d) = \MExpect\left(\#\{t \in H_{a,1}^n: d_{H_{a,1}^n}(t) = d\}\right).
$$
Let $Y_n(k)$ denote the number of $k$-vertices in $H_{a,1}^n$.

In this paper we study second degrees of vertices in $H_{a,1}^n$. The main results are stated in Theorems \ref{exp} and \ref{th3}.

We also consider the variable $X_n(k)$ equal to the number of vertices with second degree $k$ in $H_{a,1}^n$.
Note that $Y_n(k) = \sum_{i\ge k} X_n(i)$.

\subsection{Expectation}\label{ss22}
\begin{Th}\label{exp}
For any $k > 1$ we have
$$
\MExpect Y_n(k) =
\frac{(a+1)\Gamma(2a+1)}{\Gamma(a+1)k^{a}} n \left(1 +  O\left(\frac{(\ln k)^{\lceil a+1\rceil }}{k}\right)
+O\left(\frac{k^{1+a}}{n}\right) \right).
$$
\end{Th}
The easy consequence of Theorem \ref{exp} is
\begin{cor}\label {cor1}
We have $\MExpect Y_n(k) = \Theta\left(\frac n{k^a}\right)$ for $k = O\left(n^{\frac{1}{1+a}}\right)$.
\end{cor}

Using the same technique as in proof of Theorem \ref{exp} we can prove the following

\begin{Th}\label{X_n}
For any $k\ge 1$ we have
$$
\MExpect X_n(k) =
\frac{(a+1)\Gamma(2a+1)n}{\Gamma(a)k^{a+1}} \left(1 +  O\left(\frac{(\ln k)^{\lceil a+1\rceil }}{k}\right) +
O\left(\frac{k^{1+a}}{n}\right) \right).
$$
\end{Th}

Again, as a consequence we get
\begin{cor}\label {cor11}
We have $\MExpect X_n(k) = \Theta\left(\frac n{k^{1+a}}\right)$ for $k = O\left(n^{\frac{1}{1+a}}\right)$.
\end{cor}

We need the following definition.
Let $N_n(l,k)$ be the number of vertices in $H_{a,1}^n$ with degree $l$, with second degree $k$, and without loops:
$$
N_n(l,k) = \#\{t \in H_{a,1}^n: d(t)=l, d_2(t) = k, tt \notin H_{a,1}^n\}.
$$
To prove Theorem \ref{exp} we need the following auxiliary theorem.

\begin{Th}\label{th4}
In $H_{a,1}^n$ we have
$$
\MExpect N_n(l,k) = c(l,k) \,(n + \theta(n,l,k)),
$$
where $|\theta(n,l,k)| < C(l+k)^{1+a}$.
The constants $c(l,k)$ are defined as follows:
\begin{eqnarray*}
c(l,0)&=& c(0,k)= 0,\\
c(1,k)&=&c(1,k-1)\frac{a+k-1}{k+3a+1}+c(k)\frac{a+k-1}{k+3a+1},\quad k>0,\\
c(l,k)&=&c(l,k-1)\frac{al+k-1}{l(1+a)+k+2a}+c(l-1,k)\frac{l-2+a}{l(1+a)+k+2a},\quad k>0,l>1.
\end{eqnarray*}
Here $c(k) = \frac {\B(k-1+a,a+2)} {\B(a,a+1)}$.
\end{Th}

To prove these theorems we shall use two lemmas.
In \cite{Gr} Grechnikov obtained the following result.

\begin{Lem}\label{Gr}
Let $k \ge 1$ be natural; then
$$
M_n^1(k) = \frac {\B(k-1+a,a+2) n} {\B(a,a+1)} + \tilde\theta(n,k),
$$
where $|\tilde \theta(n,k)| < \tilde C/k$.
\end{Lem}

Denote by $P_n(l,k)$ the number of vertices in $H_{a,1}^n$ with a loop, with degree $l$, and with second degree $k$.

\begin{Lem}\label{lem2}
For any $n$ we have
$$
\MExpect P_n(l,k) \le p(l,k),
$$
where
\begin{eqnarray*}
p(2,0)&=&P,\\
p(l,0)&=&p(l-1,0)\frac{l-2+a}{l(1+a)-2-a},\quad l\ge3, \\
p(l,k)&=&p(l,k-1)\frac{al+k-2a-1}{l(1+a)+k-1-a}+p(l-1,k)\frac{l-2+a}{l(1+a)+k-1-a},\quad l\ge3,k\ge1.
\end{eqnarray*}
Here $P$ is some constant.
For the other values of $l$ and $k$ we have $p(l,k) = 0$.
\end{Lem}

\subsection{Concentration}
\begin{Th}\label{th3}
Let $\delta>0$ and $k=O\left(n^{\frac1{2+a+\delta}}\right)$.
Then for some $\epsilon>0$  we have
$$\Prob\left(|Y_n(k)-\MExpect(Y_n(k))|>\bigl(\MExpect(Y_n(k))\bigr)^{1-\epsilon}\right)=\bar{o}(1).$$
\end{Th}
It is a concentration result which means that the distribution of second degrees does, as the distribution of degrees, obey
(asymptotically) a power law.

This theorem is a non-trivial application of Talagrand's inequality (see \cite{Tal}). Instead of Talagrand's inequality it is possible to apply Azuma's inequality (see \cite{AlonSpencer}), but (as we show later) the result would have been weaker with Azuma's inequality.

We can prove an analogous result for the value $X_n(k)$.

\begin{Th}\label{th31}
Let $\delta>0$ and $k=O\left(n^{\frac1{4+a+\delta}}\right)$.
Then for some $\epsilon>0$  we have
$$\Prob\left(|X_n(k)-\MExpect(X_n(k))|>\bigl(\MExpect(X_n(k))\bigr)^{1-\epsilon}\right)=\bar{o}(1).$$
\end{Th}

If we substitute $a=1$ in the Buckley--Osthus model then we obtain the Bollob\'as--Riordan model \cite{BR1}.
The second degrees in this model were considered in \cite{O-G}.
The concentration of second degrees in \cite{O-G} was proved using Azuma's inequality.
This inequality provided the concentration of $X_n(k)$ around its expectation for all
$k = O\left(n^{\frac1{6+\delta}}\right)$ (with any positive $\delta$).
As stated in Theorem \ref{th31} Talagrand's inequality gives the stronger result:
for Bollob\'as--Riordan model we obtain the concentration for all $k = O\left(n^{\frac1{5+\delta}}\right)$.
We obtain this improvement in spite of the fact that the proof of the concentration of $X_n(k)$ in Theorem \ref{th31} uses the concentration of $Y_n(k)$ from Theorem \ref{th3},
so it is not optimal in this sense.

It is possible to generalize Theorem \ref{th3} (and also Theorem \ref{th31}) to the case of arbitrary  $m>1$.
The only problem in this case is that we could not prove an analog of Theorem \ref{exp} (or Corollary \ref{cor1}) for $m>1$
since it demands even more calculations. But one would expect that the following conjecture is true.

\begin{cj}\label{cj1} For any $m>1$ and $k =  O\left(n^{\min\left\{\frac1{2+a},\frac1{2a}\right\}}\right)$ we have $\MExpect Y^m_n(k) = \Theta\left(\frac n{k^a}\right)$, where  $Y^m_n(k)$ is the number of $k$-vertices in $H_{a,m}^n$.
\end{cj}

We can generalize Theorem \ref{th3} in the following way.

\begin{Th}\label{th3'} Suppose Conjecture \ref{cj1} is true. Let $m\in \mathbb{N},\ \delta>0$ and $k=O\left(n^{\min\left\{\frac1{2+a+\delta},\frac1{2a+\delta}\right\}}\right)$.
Then for some $\epsilon>0$  we have
$$\Prob\left(|Y^m_n(k)-\MExpect(Y^m_n(k))|>\bigl(\MExpect(Y^m_n(k))\bigr)^{1-\epsilon}\right)=\bar{o}(1).$$
\end{Th}

In Subsections \ref{ss31} -- \ref{ss34} we prove Theorem \ref{th3} (using  Corollary \ref{cor1}). In Subsection \ref{ss341} we prove Theorem \ref{th31} using Corollary \ref{cor11}.
In Subsection \ref{ss35} we present the sketch of the proof of Theorem \ref{th3'}.
In Section \ref{sec4} we prove results from Subsection \ref{ss22} (Theorem \ref{exp}, Theorem \ref{th4},  and Lemma \ref{lem2} in Subsections \ref{ss42}, \ref{ss41} and \ref{ss43} respectively). Finally, we prove Theorem \ref{X_n} in Subsection \ref{Corol}.

\section{Concentration}\label{sec3}

\subsection{Interpretation of the Buckley--Osthus model \\ in terms of independent variables}\label{ss31}
We consider the following sequence:
$$1, \xi_1,2,\xi_2,\ldots, n,\xi_n,$$
where $\xi_1, \dots, \xi_n$ are mutually independent random variables.
For every $ i $, we have $\xi_i: \Omega_i \rightarrow \{1,\ldots,2i-1\}$
(here $(\Omega_i, \mathcal{F}_i, \Prob_i)$ is some probability space) and
$$\Prob_i(\xi_i=2j-1)=\frac a {(a+1)i-1} \,\,\,\,\, \forall j=1,\ldots, i,$$
$$\Prob_i(\xi_i=2j)=\frac 1 {(a+1)i-1} \,\,\,\,\, \forall j=1,\ldots, i-1.$$

We can interpret  the sequence in the following way. Each $i$ is a vertex of a graph. Each  $\xi_i$ is an endpoint of the edge that goes from the vertex $i.$ If $\xi_i=2j-1,$ then the edge goes to the vertex $j$. If $\xi_i=2j$, then we say that the edge from the vertex $i$ goes to the same vertex as the edge from the vertex $j$. The value of the variable $\xi_j$ can also be even (say $\xi_j=2j_1$, for some integer $j_1$), then the edge from the vertex $i$ is again redirected according to the variable
$\xi_{j_1}$. Finally this process stops at some odd value $2v-1$ and we say that $\xi_i$ (as well as $\xi_j$ and $\xi_{j_1}$) {\it leads} to the vertex $v$. We also say that $\xi_i$ leads to $\xi_j.$

It is not hard to check that the graph model we obtained is exactly the Buckley--Osthus model.
Indeed, at each time step $i$
the in-degree of each vertex $ j \in \{1, \dots, i-1\} $ is equal to the number of  variables that lead (directly or indirectly) to the vertex $j$.

Let us give yet another interpretation of the model described above.
Consider a vertex $v$ from the obtained graph.
We can think of
all the variables that lead to $v$ as connected as a rooted tree, with $v$ as the root.
Let $X = \{\xi_{i_1}, \ldots, \xi_{i_d}\}$ be the set of variables that lead to $v$.
We inductively construct the corresponding tree on $d$ vertices $i_1,\ldots, i_d$.
First consider those  variables $\xi_{i_1^1},\ldots, \xi_{i_{l_1}^1}$ from $X$ that lead to $v$ directly.
The corresponding vertices $i_1^1, \ldots, i_{l_1}^1$
are adjacent to $v$ in the tree.
Suppose we choose all the vertices at distance $\le s$ from $v$
and $i_1^s,\ldots, i_{l_s}^s$ are the vertices at distance $s$ from $v$.
Consider the set $\{\xi_{i_1^{s+1}},\ldots, \xi_{i_{l_{s+1}}^{s+1}}\}$ of variables that lead to some of
$\xi_{i_1^s},\ldots, \xi_{i_{l_s}^s}.$
We join each of the vertices $i_1^{s+1},\ldots, i_{l_{s+1}}^{s+1}$ to
the corresponding vertex  from $\{i_1^s,\ldots, i_{l_s}^s\}.$
We thus obtain
the set of vertices at distance $s+1$ from $v$.

\subsection{Decreasing the number of $k$-vertices}

We fix a value $\mathbf{x}=(x_1,\ldots,x_n)$
of the random vector $\mathbf{\xi}=(\xi_1,\ldots, \xi_n)$ from the probability space $\Omega=\prod_{i=1}^n\Omega_i$.
The quantity $Y_n(k)$ is
a function from $\Omega$ to $\mathbb{N}.$  We discuss the following question. How
can the value $Y_n(k)=Y_n(k,\mathbf{x})$ decrease, if we change one coordinate $x_i$ of the vector $\mathbf{x}$?
In other words, we want to find $c(i,\mathbf{x}) = \max_{{\bf x}'} (Y_n(k,\mathbf{x})-Y_n(k,\mathbf{x'})),$ where
$\mathbf{x'}$ is an arbitrary vector that differs from $\mathbf{x}$ in exactly the $i$th coordinate.

\begin{Lem}\label{lem1}
For any $\mathbf{x}=(x_1,\ldots, x_n)$ and $ i\in \{1, \dots, n\}$ we have $c(i,\mathbf{x})\le 2k+1.$
\end{Lem}

\begin{Proof}
It is convenient to think about the tree interpretation of the random variables. If we change a value $x_i$ of one random variable
$\xi_i$ to some value $x'_i$, then
all the variables that lead to $x_i$ are redirected to $x'_i$.
In terms of the tree interpretation, we pick the branch of the tree in which all the edges lead to $x_i$:
if $x'_i$ is odd, then we link the branch to the vertex with the number $(x'_i+1)/2$;
if $x'_i$ is even,  then we link the branch to the variable 
$\xi_{x'_i/2}$.

We want to interpret the change of one coordinate in terms of the graph $H_{a,1}^{n}.$
Suppose $x_i$ leads to a vertex $v$. Then all the variables that lead to $x_i$ lead to $v$.
If we change $x_i$ to $x_i'$ and $x_i'$ leads to $v'$, then we change the value of all such variables from $v$ to $v'$.
Or, in terms of $H_{a,1}^{n}$, we take a bundle of edges in the vertex $v$ and move the bundle to the vertex $v'$.
More precisely, if we had a bundle of edges $(i_1,v), \ldots, (i_d,v)$, then after the change we have the edges
$(i_1,v'), \ldots, (i_d,v')$. All the rest stays the same.

Now we go on to the proof. We should show that after
the change of the $i$th coordinate,
the number of $k$-vertices
we spoil does not exceed $2k+1$.
Suppose we moved a bundle of edges $(i_1,v), \ldots, (i_d,v)$.
It is easy to see that we could spoil only the $k$-vertices that have a common edge with $v$ or $v$ itself.
Note that we could not spoil the $k$-vertices in the neighborhood of $v'$.

We split the set $N_v$ of the vertices incident to $v$ into two parts: $I=\bigcup_{j=1}^d \{i_j\}$ and $N_v\backslash I$.
If $|N_v\backslash I|\ge k+1,$ then after changing the edges from the bundle, all the $k$-vertices from $N_v\backslash I$ are still $k$-vertices.
Indeed, all the edges in vertex $v$ except for one are $2$-incident edges for any neighbor of $v$, so there are at least $k$ such edges for every
vertex from $N_v\backslash I$.
Similarly, if $|I|\ge k+1$, then no $k$-vertices among $i_1,\ldots, i_d$ are spoiled except for at most one,
since they are all adjacent to the vertex $v'$.
The only case when some of $i_1,\ldots, i_d$ is spoiled is $i_j = v'$
and so we will not count the edges $(i_1,v'), \ldots, (i_d,v')$ in the second degree of $i_j$.

Finally, the number of $k$-vertices we spoil
does not exceed $\min\{|I|,k\}+\min\{|N_v\backslash I|,k\}+1\le 2k+1$.

\end{Proof}

We now want to estimate the influence of each variable more accurately.
Suppose $Y_n(k,\mathbf{x})=q.$
For each $k$-vertex $v_j$, $j=1,\ldots, q,$ we consider a subset of coordinates $K_j=K_{v_j}(\mathbf{x}) = \{i_1^j, \dots, i_{d_j}^j\}$,
such that $v_j$ is a $k$-vertex for any $\mathbf{y}$ that agrees with $\mathbf{x}$ on the coordinates from $K_j$. It is worth noting that
$ K_j $ depends on $ {\bf x} $, but is not uniquely defined by it. For any choice of the sets $ K_1, \dots, K_q $,
we denote their collection by $\mathcal{K}=\mathcal{K}(\mathbf{x})$.
Clearly, $Y_n(k,\mathbf{y})\ge q$, for any $\mathbf{y}$ that agrees with $\mathbf{x}$ on all the coordinates from all $K_j \in \mathcal{K}$.



For each coordinate $i$, we define its multiplicity $C(i,\mathbf{x},\mathcal{K})= |\{j:i\in K_j\}|.$

 It is easy to see that for any $ {\bf x} $ and any $\mathcal{K}$ one has $c(i,\mathbf{x})\le C(i,\mathbf{x},\mathcal{K}).$
 So we have $$c(i,\mathbf{x})\le \min\{2k+1, C(i,\mathbf{x},\mathcal{K})\}=:c_i(\mathbf{x},\mathcal{K}).$$

We call a collection $\mathcal{K}$ \textit{stable}, if for every $k$-vertex $v_i$ we construct all the sets $K_i$ according to the following rule:
if $K_i$ contains some of the coordinates that lead to a vertex $w$, then $K_i$ contains all of the coordinates that lead to $w$.

Now, it should not be surprising that for such set systems we can prove an analog of Lemma \ref{lem1}. Namely, consider a vector
$\mathbf{x}$, the corresponding $k$-vertices $v_j,j=1,\ldots, q,$ and some stable collection $\mathcal{K}$.
Let $i\in \{1, \dots, n\}$ and
$\mathcal{K} \setminus\{i\}:=\{K_j\setminus\{i\},j=1,\dots,q\}.$ Given $i$
there exist at least
$q-c_i(\mathbf{x},K)$ such $k$-vertices that are $k$-vertices for any $\mathbf{x}'$ with $x'_s=x_s$ for all
$s\in K_j\setminus\{i\},j=1,\dots,q$.
To prove this fact one has to follow the proof of Lemma \ref{lem1} and make sure that the proof works also for this case.
The only additional consideration needed is that the number of $k$-vertices we loose does not exceed the multiplicity $C(j,\mathbf{x},\mathcal{K})$.

\begin{Lem}\label{lemadd}
Let $\mathcal{K}$ be a stable collection as described above.
We have $Y_n(k,\mathbf{x})-Y_n(k,\mathbf{x'})\le \sum_{j\in J}c_j(\mathbf{x},\mathcal{K})$
for any vector $\mathbf{x}'$ such that $x'_i=x_i$ for all $i\in \{1, \dots, n\}\backslash J.$
\end{Lem}

 \begin{Proof}
Suppose we change one coordinate $j$ of $\mathbf{x}$ and obtain some vector $\mathbf{\hat{x}}$.
Then we consider $d_j := C(j,\mathbf{x},\mathcal{K})$ $k$-vertices $w_1,\ldots, w_{d_j}$
such that $j\in K_{w_i}(\mathbf{x}).$ We remove $j$ from each of these sets and
we check for each $i=1,\ldots, d_j$ whether the obtained collection guarantees $w_i$ to be a $k$-vertex or not.
If $w_i$ is a $k$-vertex then we define $K_{w_i}(\mathbf{\hat{x}}) = K_{w_i}(\mathbf{x})\backslash\{j\}.$
If $w_i$ is not a $k$-vertex then we exclude the set $K_{w_i}(\mathbf{x})$ from $\mathcal{K}.$
At the end of this step we obtain a new collection $\mathcal{\hat{K}}.$
To prove the lemma we need one consideration. Namely, instead of changing the edges $(i_1,v), \ldots, (i_d,v)$ to $(i_1,v'), \ldots, (i_d,v')$
we can create a new imaginary vertex $w$ and change the edges to $(i_1,w), \ldots, (i_d,w).$
We denote the obtained graph by $G_w$. We do not count $w$ as a $k$-vertex even if it has $\ge k$ 2-incident edges.
It is easy to check that for this graph the collection $\mathcal{\hat{K}}$ is stable.

The number of $k$-vertices (except for $w$) in the graph $G_w$ is definitely not bigger than the same number for the graph
corresponding to $\mathbf{\hat{x}}$. Moreover, the multiplicity of each coordinate in $\mathcal{\hat{K}}$ is less than or equal to
the corresponding multiplicity in $\mathcal{K}$. We also have $Y_n(k,\mathbf{x})-Y_n(k,G_w)\le c_j(\mathbf{x},\mathcal{K})$.
Similarly, if $\mathbf{x'}$ differs from $\mathbf{x}$ in $l$ coordinates, then the graph corresponding to $\mathbf{x'}$ has at least
as many $k$-vertices as the graph $G$ obtained by forming $l$ imaginary vertices. Moreover, at each step (if we change the coordinate $j'$
and form the corresponding graph $G'$) we spoil at most $\min\{2k+1, C(j',\mathbf{x},\mathcal{K})\}$ $k$-vertices and obtain a stable set system.

Consequently, we have $Y_n(k,\mathbf{x})-Y_n(k,\mathbf{x'})\le Y_n(k,\mathbf{x})-Y_n(k,G)\le \sum_{j\in J}c_j(\mathbf{x},\mathcal{K}).$

\end{Proof}

\subsection {Construction of a suitable set $\mathcal{K}$}

\begin{Lem}\label{lemK} Suppose $Y_n(k,\mathbf{x})=q$ for some vector $\mathbf{x}$ in the corresponding graph $G_{\bf x}$.
Then we can construct a stable set system $\mathcal{K}=\{K_1,\ldots, K_q\}$ such that $\sum_{i=1}^n c_i(\mathbf{x},\mathcal{K})\le (4k+5)q.$
\end{Lem}

\begin{Proof} First consider the set $V$ of vertices with degree at least $k+2$.
Put $NV=\{u: u$ is a neighbor of $v\in V\}$.
Note that a vertex from $V$ can also belong to $NV$.
Assume that $|NV| = z$.
All vertices from $NV$ are $k$-vertices.
Let $BV$ be the set of vertices from $NV$ which do not have an outcoming edge that goes to $V$.
We have $|BV|\le z/(k+1)$ since each vertex has at most one outcoming edge.

We denote by $L_v, L_v \subset \{1, \dots, n\}$, the set of coordinates that lead to $v$.
We also put $LV = \cup_{v\in V}L_v$ and $LBV = \cup_{v\in BV}L_v$.
For any $u \in NV$ we put $K_u=LV \cup LBV$.
It is easy to see that for any $\mathbf{x'}$ such that $x_i = x'_i$ for every $i\in LV\cup LBV$, the vertex $u$ is
$k$-vertex in the graph corresponding to $\mathbf{x'}.$

For $i \in LV \cup LBV$ we estimate $c_i(\mathbf{x},\mathcal{K})$ by $2k+1$.
Note that $|LV| \le z + \frac{z}{k}$.
We add additional $\frac z k$ variables because the vertices from $V$ can have loops.
We have $\deg w \le k+1$ for $w\in BV\backslash V$ and $|BV\backslash V|\le z/(k+1)$, therefore $|LBV\backslash LV|\le z$.
So we have
$$
\sum_{i\in LV\cup LBV} c_i(\mathbf{x},\mathcal{K})\le (2k+1)(|LBV\backslash LV|+|LV|)\le (2k+1)\left(2z+\frac zk\right) \le (4k+5)z.
$$


Next we consider the set $W$ of the remaining  $k$-vertices.
We have $|W|=q-z$.
By the definition, for  any $w\in W$ all the neighbors $N_w$ of $w$ have degree less than or equal to $k+1$.

For each $w\in W$ we consider $V_w=\{v_1, \ldots, v_w\},$ $V_w\subset N_w,$
such that the number of edges adjacent to at least one of $v_i \in N_w$ and not adjacent to $w$
is between $k$ and $2k$. We can find such $V_w$ since $w$ is a $k$-vertex. We can choose $v_i\in N_w$ one by one,
until the total number of $2$-adjacent edges does not exceed $k$. But it cannot exceed $2k$ since  $\deg v_i\le k+1$
for $ v_i\in N_w$. Denote by $LV_w$ all the variables that lead to $V_w$.

Now for each $w\in W$ we put $K_w = LV_w\cup L_w$.
Note that $|LV_w|\le 2k$ and for $w \in W\backslash V$ we have  $|L_w|\le k+1$.

Now we can make the final estimate:

\begin{multline}
\sum_{i=1}^n c_i(\mathbf{x},\mathcal{K})=\sum_{i\in LV\cup LBV} c_i(\mathbf{x},\mathcal{K})
+ \sum_{i\in \{1, \dots, n\}\backslash {LV \cup LBV}} c_i(\mathbf{x},\mathcal{K})
\le (4k+5)z+\\ + \sum_{w \in W} |LV_w|+\sum_{w \in W\backslash V} |L_w|\le (4k+5)z+ 2k(q-z)+(k+1)(q-z)\le (4k+5)q.\notag
\end{multline}

The fact that $\mathcal{K}$ is a stable set system follows from the construction.
\end{Proof}

\subsection {Application of Talagrand's inequality}\label{ss34}
First we briefly review  Talagrand's inequality (see e.g., \cite{AlonSpencer}).

Let $\Omega = \prod_{i=1}^n\Omega_i$ be a product probability space with product measure. Suppose $\mathbf{\alpha} = (\alpha_1,\ldots, \alpha_n)$, $\sum_{i=1}^n \alpha_i^2=1.$
We define the following distance between a set $A\subset \Omega$ and a point $\mathbf{x}\in \Omega$:
$$\dist (A,\mathbf{x})=\max_{\mathbf{\alpha}}\min_{\mathbf{y}\in A}\sum_{i\in I_{{\bf x}{\bf y}}}\alpha_i,$$
where $I_{{\bf x}{\bf y}} = \{i: x_i\ne y_i\}.$

For $t>0$ we denote by $A_t$ the set $\{\mathbf{x}: \dist (A,\mathbf{x})\le t\}.$

\begin{Th}{\bf(Talagrand's inequality)} For any $t>0$ and $A\subset \Omega$ we have $\Prob(A)(1- \Prob(A_t))\le e^{-\frac {t^2}4}.$
\end{Th}

We use the inequality to derive the following theorem.
\begin{Th}\label{tal}
For $t>0,$ $k, s\in \mathbb{N}$, and $f(s)$ that satisfies the condition  $f^2(s)> (2k+1)(4k+5) s$ we have $\Prob\bigl(Y_n(k)\le s-t f(s)\bigr) \Prob\bigl(Y_n(k)\ge s\bigr)\le e^{-\frac {t^2}4}$.
\end{Th}

\begin{Proof}
The inequality is trivial for $tf(s)> s,$ so we can assume w.l.o.g. that $tf(s)\le s$.  Since $\xi_i$ are independent,
we can apply Talagrand's inequality to the points $\mathbf{x}$ from the probability space $\Omega$.
Actually, all we need to prove is that for any $\mathbf{x}$ such that $ Y_n(k,\mathbf{x})\ge s$ we have ${\bf x}\notin A_t$,
where $A = \{\mathbf{y}: Y_n(k,\mathbf{y})\le s -t f(s) \}$.

Suppose $Y_n(k,\mathbf{x})= q\ge s.$ Given $\mathbf{x}$ we fix a set system $\mathcal{K}$ as in Lemma \ref{lemK}.
Then by Lemma \ref{lemadd} for any ${\bf y}\in A$ we have $q-s+ t f(s)\le Y_n(k,\mathbf{x})-Y_n(k,\mathbf{y})\le \sum_{j\in I_{{\bf x}{\bf y}}}c_j(\mathbf{x},\mathcal{K})$.

We define a suitable vector $\mathbf{\alpha}=\mathbf{\alpha}(\mathbf{x})$. Namely, $\alpha_i = \frac{c_i(\mathbf{x},\mathcal{K})}{\sqrt{\sum_{j=1}^{n} c_j^2(\mathbf{x},\mathcal{K})}}.$ It is easy to see that $\sum\alpha_j^2 = 1.$

We have
$$\sum_{j=1}^{n} c_j^2(\mathbf{x},\mathcal{K})\le \max_j c_j(\mathbf{x},\mathcal{K})\sum_{j=1}^{n} c_j(\mathbf{x},\mathcal{K})\le (2k+1)(4k+5)q.$$
In the last inequality we used Lemma \ref{lemK} and the definition of $c_j(\mathbf{x},\mathcal{K})$.

Now we show that  $\sum_{i\in I_{{\bf x}{\bf y}}}\alpha_i >t$ for any ${\bf y}\in A$.
\begin{equation}\label{est1}\sum_{i\in I_{{\bf x}{\bf y}}}\alpha_i =\frac{\sum_{i\in I_{{\bf x}{\bf y}}}c_i(\mathbf{x},\mathcal{K})}{\sqrt{\sum_{j=1}^{n} c_j^2(\mathbf{x},\mathcal{K})}}\ge
\frac{q-s+ t f(s)}{\sqrt{(2k+1)(4k+5)q}}\ge \frac{t f(s)}{\sqrt{(2k+1)(4k+5)s}}>t.
\end{equation}
The second inequality holds since for $q\ge s, tf(s)\le s$ we have $\frac {q-s+tf(s)}q\ge \frac{tf(s)}s.$
The last inequality follows from the statement of the theorem.

From \eqref{est1} we obtain that $\dist (A,\mathbf{x})>t$, in other words, ${\bf x}\notin A_t.$
\end{Proof}

 We apply Theorem \ref{tal} with $t=2\ln n$, $s= m(Y_n(k)) +t (\MExpect Y_n(k))^{1-\varepsilon}$, $f(s) = (\MExpect Y_n(k))^{1-\varepsilon}.$
Here $m(Y_n(k))$ is the median of $Y_n(k),$ and, consequently, $\Prob\bigl(Y_n(k)\le s-t f(s)\bigr)\ge 1/2.$
Since for any random variable $Z$ we have $m(Z)\le 2\MExpect Z,$  it is easy to see that the conditions of
Theorem \ref{tal} hold if $$(\MExpect Y_n(k))^{1-2\varepsilon}\ge 12(2k+1)^2 \ln n.$$
If $\varepsilon$ is small enough then this inequality is a consequence of Corollary \ref{cor1} and the conditions of Theorem \ref{th3}.

 We obtain that $$\Prob\Bigl(Y_n(k)\ge m(Y_n(k)) +2\ln n (\MExpect Y_n(k))^{1-\varepsilon}\Bigr)\le 2 e^{-\frac {t^2}4} = \bar{o}(1/n),$$
 and since $Y_n(k)\le n$ for all $k,$ we have  $$\MExpect Y_n(k)\le m(Y_n(k)) +2\ln n (\MExpect Y_n(k))^{1-\varepsilon} + \bar{o}(1).$$

 Similarly we can derive that $$\Prob\Bigl(Y_n(k)\le m(Y_n(k)) - 2\ln n (\MExpect Y_n(k))^{1-\varepsilon}\Bigr)\le
2 e^{-\frac {t^2}4} = \bar{o}(1/n)$$ and $$\MExpect Y_n(k)\ge m(Y_n(k)) - 2\ln n (\MExpect Y_n(k))^{1-\varepsilon} - \bar{o}(1).$$
Consequently, for some $\delta >0$ and all sufficiently large $n$ we have $|\MExpect Y_n(k)- m(Y_n(k))| \le (\MExpect Y_n(k))^{1-\delta}$.
Therefore, for some $\varepsilon'>0$

$$\Prob\left(|Y_n(k)-\MExpect(Y_n(k))|>\MExpect(Y_n(k))^{1-\varepsilon'}\right)=\bar{o}(1).$$
This concludes the proof of Theorem \ref{th3}.

\subsection {Proof of Theorem \ref{th31}}\label{ss341}
We use the obvious fact that $X_n(k)=Y_n(k)-Y_n(k+1).$ Fix some $\epsilon'>0.$
First we want to apply Theorem \ref{tal} to $Y_n(k)$ and $Y_n(k+1).$ We argue as after the proof of Theorem \ref{tal}. We put $f(s)=\frac{n^{1-\epsilon'}}{k^{1+a}},$ $t=2\ln n$ and  $s_1= m(Y_n(k)) +tf(s), s_2=m(Y_n(k)), s_3=m(Y_n(k+1))+tf(s), s_4=m(Y_n(k+1)).$

We apply Theorem \ref{tal} to $Y_n(k)$ with $s_1$ and $s_2$, and to $Y_n(k+1)$ with $s_3$ and $s_4$ and obtain
$$\Prob\left(|Y_n(k)-\MExpect(Y_n(k))|>\frac{n^{1-\epsilon'+\bar{o}(1)}}{k^{1+a}}\right)=\bar{o}(1),$$
$$\Prob\left(|Y_n(k+1)-\MExpect(Y_n(k+1))|>\frac{n^{1-\epsilon'+\bar{o}(1)}}{k^{1+a}}\right)=\bar{o}(1),$$
provided $$\frac{n^{2-2\epsilon'}}{k^{2+2a}}\ge \Theta\left(nk^{2-a}\right) \ln n.$$

It is easy to see that this holds if the conditions of Theorem \ref{th31} are satisfied for some $\delta>0.$ We have $|X_n(k)-\MExpect(X_n(k))|\le |Y_n(k)-\MExpect(Y_n(k))|+|Y_n(k+1)-\MExpect(Y_n(k+1))|,$ so
$$\Prob\left(|X_n(k)-\MExpect(X_n(k))|>\frac{n^{1-\epsilon'+\bar{o}(1)}}{k^{1+a}}\right)=\bar{o}(1).$$

Since $\frac{n^{1-\epsilon'+\bar{o}(1)}}{k^{1+a}}=\MExpect(X_n(k))^{1-\epsilon}$ for some $\epsilon>0$, this inequality completes the proof of Theorem \ref{th31}.

\subsection{Generalization to the case of arbitrary $m$}\label{ss35}
The proof of Theorem \ref{tal} can be modified to the case of the graph $H_{a,m}^n$. We present only the sketch of the argument. Suppose $m>1$ is fixed. The number of variables changes from $n$ to $mn$. The interpretation in terms of independent variables works for this case. Lemmas \ref{lem1}, \ref{lemadd}, \ref{lemK} hold for $m>1$, but with some minor changes.

When we take a bundle of edges from a vertex $v$ and move it to some vertex $v'$, we can spoil not only the neighborhood of $v$, but also the vertex $v'$. Namely, suppose we change edges $(v,w_1),\ldots, (v,w_l)$ to $(v',w_1),\ldots, (v',w_l)$. If $v'$ was a $k$-vertex and edges $(v,w_1),\ldots, (v,w_l)$ were counted in the second degree of $v'$, then the second degree of $v'$ may decrease (this is impossible in the graph $H_{a,1}^n$ since $H_{a,1}^n$ is a tree). It is not difficult to see that we cannot spoil the other vertices.

 Hence, we can formulate some analogs of  Lemmas \ref{lem1}, \ref{lemadd}.
\begin{Lem}\label{lem1'}
For any $\mathbf{x}=(x_1,\ldots, x_n)$ and $ i\in \{1, \dots, n\}$ we have $c(i,\mathbf{x})\le 2k+2.$
\end{Lem}

 We also have $$c(i,\mathbf{x})-1\le \min\{2k+1, C(i,\mathbf{x},\mathcal{K})\}=:c_i(\mathbf{x},\mathcal{K}).$$

\begin{Lem}\label{lemadd'}
 Let $\mathcal{K}$ be a stable set collection. We have $Y_n(k,\mathbf{x})-Y_n(k,\mathbf{x'})\le \sum_{j\in J}(c_j(\mathbf{x},\mathcal{K})+1)$ for any vector $\mathbf{x}'$ such that $x'_i=x_i$ for all $i\in \{1, \dots, n\}\backslash J.$
\end{Lem}

  Lemma \ref{lemK} holds for $c_j(\mathbf{x},\mathcal{K})$ and $m>1$ without any changes.

The only thing left is to modify the proof of Theorem \ref{tal}. We put $\alpha_i = \frac{c_i(\mathbf{x},\mathcal{K})+1}{\sqrt{\sum_{j=1}^{mn} (c_j(\mathbf{x},\mathcal{K})+1)^2}}.$ Then

$$\sum_{j=1}^{mn} (c_j(\mathbf{x},\mathcal{K})+1)^2\le \left(\max_j c_j(\mathbf{x},\mathcal{K})+2\right)\sum_{j=1}^{mn}
c_j(\mathbf{x},\mathcal{K})+mn\le (2k+3)(4k+5)q+mn.$$
Finally,
\begin{equation*}
\sum_{i\in I_{{\bf x}{\bf y}}}\alpha_i =\frac{\sum_{i\in I_{{\bf x}{\bf y}}}(c_i(\mathbf{x},\mathcal{K})+1)}{\sqrt{\sum_{j=1}^{mn}
(c_j(\mathbf{x},\mathcal{K})+1)^2}}\ge \frac{q-s+ t f(s)}{\sqrt{(2k+3)(4k+5)q+mn}}\ge \frac{t f(s)}{\sqrt{(2k+3)(4k+5)s+mn}}>t,
\end{equation*}
if $f^2(s)>(2k+3)(4k+5)s+mn.$ So we can formulate an analog of Theorem \ref{tal}.

\begin{Th}\label{tal2}
For $t>0,$ $m, k, s\in \mathbb{N}$, and $f(s)$ that satisfies the condition  $f^2(s)> (2k+3)(4k+5) s +mn$ we have $$\Prob\bigl(Y^m_n(k)\le s-t f(s)\bigr) \Prob\bigl(Y^m_n(k)\ge s\bigr)\le e^{-\frac {t^2}4}.$$
\end{Th}

Finally, arguing as after the proof of Theorem \ref{tal}, one can see that Theorem \ref{th3'} follows from Theorem \ref{tal2} and Conjecture \ref{cj1}.
\section{Estimation of $\MExpect Y_n(k)$}\label{sec4}

We need the following
notation.
Let $X$ be a function on $n$ (the number of vertices), $l$ (the first degree we are interested in),
$k$ (the second degree we are interested in);
then denote by $\theta(X)$ some function on $n$, $l$, $k$ such that
$|\theta(X)| < X$.

\subsection{Proof of Theorem \ref{th4}}\label{ss41}

It follows from the definition of $H_{a,1}^n$ that $N_n(l,0) = N_n(0,k)=0$.
Indeed, since we have no vertices of degree $0$, we see that $N_n(0,k)=0$.
Since vertices with loops are not counted in $N_n(l,k)$,
 we have no vertices of second degree $0$ and $N_n(l,0) =0$.
Therefore, we have $\MExpect N_n(l,0)=\MExpect N_n(0,k)=0$.
We want to prove that there exists such constant $C$ that
$$
\MExpect N_n(l,k) = c(l,k) \,(n + \theta(n,l,k)),
$$
where $|\theta(n,l,k)| < C(l+k)^{1+a}$.

Let us demonstrate that
$\MExpect N_n(1,k) = c(1,k) \,\left(n + \theta\left(C k^{1+a}\right)\right)$.
We shall use induction on $k$.
For $k=0$ there is nothing to prove.

Assume that for $j<k$ we have
$$
\MExpect N_n(1,j) = c(1,j) \,(n + \theta\left(C j^{1+a}\right)).
$$

Denote by $N_i(l)$ the number of vertices with degree $l$ in $H_{a,1}^i$.
We use induction on $i$ and the equality
\begin{multline}\label{l=1}
\MExpect ( N_{i+1} (1,k) | N_{i}(1,k),N_{i}(1,k-1),N_{i}(k)) =\\
= N_i(1,k) \left( 1 - \frac {k+2a} {(a+1)i + a} \right) + \frac {(k - 1 + a) N_i(1,k-1)} {(a+1)i + a}
+ \frac {(k - 1 + a) \, N_i(k)} {(a+1)i + a}.
\end{multline}
Let us explain this equality.
Suppose we constructed $H_{a,1}^i$.
We add one vertex and one edge.
There are $N_i(1,k)$ vertices with degree $1$ and with second degree $k$ in $H_{a,1}^i$.
The probability that we ``spoil'' one of these vertices is $\frac{k+2a}{(a+1)i + a}$.
We also have $N_i(1,k-1)$ vertices with degree $1$ and with second degree $k-1$.
The probability that one of these vertices has degree $1$ and second degree $k$ in $H_{a,1}^{i+1}$ is
$\frac{k - 1 + a}{(a+1)i + a}$.
Finally, with probability equal to $\frac{k - 1 + a}{(a+1)i + a}$ the vertex $i+1$ has necessary degrees in $H_{a,1}^{i+1}$.
From \eqref{l=1} we obtain
\begin{equation}\label{M1}
\MExpect N_{i+1} (1,k) =
\MExpect N_i(1,k) \left( 1 - \frac {k+2a} {(a+1)i + a} \right) + \frac {(k - 1 + a)\MExpect N_i(1,k-1)} {(a+1)i + a}
+ \frac {(k - 1 + a) \,\MExpect N_i(k)} {(a+1)i + a}.
\end{equation}

Note that if we have at least one vertex with first degree $1$ and second degree $k$ in $H_{a,1}^i$,
then we have at least $k$ edges in this graph.
Therefore $\MExpect N_{i}(1,k) = 0$ when $i < k$.
Consider the case $i = k$.
First, note that
$$
\MExpect N_{k} (1,k) \ge c(1,k)  \left(k + \theta\left(C k^{1+a}\right)\right)
$$
with some $C$.
For a finite number of small $k$ we can find a constant $C$ such that
$$
\MExpect N_{k} (1,k) = c(1,k)  \left(k + \theta\left(C k^{1+a}\right)\right).
$$

Using \eqref{M1}, Lemma \ref{Gr}, and the assumptions of the theorem we get
\begin{gather*}
\MExpect N_{k} (1,k) =
\MExpect N_{k-1}(1,k-1) \frac {k - 1 + a} {ak + k - 1}
+ M_{k-1}^1(k) \frac {k - 1 + a} {ak + k - 1} =\\
= c(1,k-1)\frac {k - 1 + a} {ak + k - 1}\left(k-1 + \theta\left(C\left(k - 1\right)^{1+a}\right)\right)
+ c(k) \frac {k - 1 + a} {ak + k - 1} \left(k-1 + \theta \left(C_1 k^{1+a}\right)\right) =\\
= c(1,k)\frac{(k+3a+1)(k-1)}{ak + k - 1} +
c(1,k-1)\frac {k - 1 + a} {ak + k - 1}\theta\left(C\left(k - 1\right)^{1+a}\right) +
c(k) \frac {k - 1 + a} {ak + k - 1} \theta \left(C_1 k^{1+a}\right) =\\
= k c(1,k) + \frac{3ak+k-3a-1-ak^2}{ak + k - 1}c(1,k) +
c(1,k-1)\frac {k - 1 + a} {ak + k - 1}\theta\left(C\left(k - 1\right)^{1+a}\right) +
\\
+c(k) \frac {k - 1 + a} {ak + k - 1} \theta \left(C_1 k^{1+a}\right) \le k c(1,k)
+ \frac{(3ak+k-3a-1-ak^2)(a+k-1)}{(ak + k - 1)(k+3a+1)}c(1,k-1)
+ \\
+ \frac{(3ak+k-3a-1-ak^2)(a+k-1)}{(ak + k - 1)(k+3a+1)}c(k) + c(1,k-1)\frac {C(k - 1 + a)} {ak + k - 1}\left(k - 1\right)^{1+a} +\\
c(k) \frac {k - 1 + a} {ak + k - 1} C_1 k^{1+a} \le k c(1,k) + \frac{C(a+k-1)}{k+3a+1}c(1,k-1)k^{1+a}
+ \frac{C(a+k-1)}{k+3a+1}c(k)k^{1+a}.
\end{gather*}
This holds for big values of $k$. Indeed,
$$
\frac{(3ak+k-3a-1-ak^2)}{(ak + k - 1)(k+3a+1)}
+\frac  {C\left(k - 1\right)^{1+a}}{ak + k - 1} \le
\frac{C}{k+3a+1}k^{1+a},
$$
if $k$ and $C$ are big enough.

Consider the case $i \ge k$.
Using \eqref{M1}, Lemma \ref{Gr}, and the inductive assumption we get
$$
\MExpect N_{i+1} (1,k) =
\MExpect N_i(1,k) \left(1 - \frac {k+2a} {(a+1)i + a} \right)
+ \MExpect N_i(1,k-1) \frac {k - 1 + a} {(a+1)i + a}
+ M_i^1(k) \frac {k - 1 + a} {(a+1)i + a} =
$$
$$
= c(1,k) \left(i + \theta \left(C k^{1+a}\right)\right) \left( 1 - \frac {k+2a} {(a+1)i + a} \right)
+ c(1,k-1) \left(i + \theta \left(C (k-1)^{1+a}\right)\right) \frac {k - 1 + a} {(a+1)i + a} +
$$
$$
+ c(k)\left(i + \theta_1 \left(C_1 k^{1+a}\right)\right) \frac{k - 1 + a} {(a+1)i + a} = c(1,k) (i + 1) -
$$
$$
- c(1,k)\frac {i(k+3a+1) + a} {(a+1)i + a}
+ c(1,k)\theta\left(C k^{1+a}\right)\left( 1 - \frac {k+2a} {(a+1)i + a} \right)
+ c(1,k-1) i  \frac {k - 1 + a} {(a+1)i + a} +
$$
$$
+ c(1,k-1)\theta \left(C (k-1)^{1+a}\right) \frac {k - 1 + a} {(a+1)i + a}
+ c(k) i \frac {k - 1 + a} {(a+1)i + a}
+ c(k) \theta_1 \left(C_1 k^{a+1}\right) \frac {k - 1 + a} {(a+1)i + a} =
$$
$$
= c(1,k) (i + 1)
+ c(1,k)\theta\left(C k^{1+a}\right)\left( 1 - \frac {k+2a} {(a+1)i + a} \right)
- \frac {(k - 1 + a)ac(1,k-1)} {((a+1)i + a)(k+3a+1)} -
$$
$$
- \frac {(k - 1 + a)ac(k)} {((a+1)i + a)(k+3a+1)}
+ c(1,k-1)\theta \left(C (k-1)^{1+a}\right) \frac {k - 1 + a} {(a+1)i + a}
+ c(k) \theta_1 \left(C_1 k^{a+1}\right) \frac {k - 1 + a} {(a+1)i + a}.
$$
We want to prove that there exists a constant $C$ such that
$$
c(1,k) C k^{1+a} \frac {k+2a} {(a+1)i + a} \ge
\frac {(k - 1 + a)ac(1,k-1)} {((a+1)i + a)(k+3a+1)}
+ \frac {(k - 1 + a)ac(k)} {((a+1)i + a)(k+3a+1)} +
$$
$$
+ c(1,k-1) C (k-1)^{1+a}\frac {k - 1 + a} {(a+1)i + a}
+ c(k) C_1 k^{1+a} \frac {k - 1 + a} {(a+1)i + a}.
$$
It is sufficient to prove that the following inequalities hold:
$$
c(1,k-1) C k^{1+a} \frac {(k+2a)(k-1+a)} {((a+1)i + a)(k+3a+1)} \ge
\frac {(k - 1 + a) a c(1,k-1)} {((a+1)i + a)(k+3a+1)}+
c(1,k-1) C (k-1)^{1+a}\frac {k - 1 + a} {(a+1)i + a}
$$
and
$$
c(k) C k^{1+a} \frac {(k+2a)(k-1+a)} {((a+1)i + a)(k+3a+1)} \ge
\frac {(k - 1 + a)ac(k)} {((a+1)i + a)(k+3a+1)}
+ c(k) C_1 k^{1+a} \frac {k - 1 + a} {(a+1)i + a}.
$$
Or
$$
C k^{1+a} (k+2a) \ge
a
+  C (k-1)^{1+a} (k+3a+1)  ,
$$
$$
 C k^{1+a} (k+2a)\ge
a + C_1 k^{1+a} (k+3a+1).
$$
Note that
$$
k^{1+a} (k+2a) - (k-1)^{1+a} (k+3a+1) =
$$
$$
= k^{1+a} (k+2a) - (k^{1+a} - (1+a)k^a + \frac{a(a+1)}{2}k^{a-1} + O(k^{a-2})) (k+3a+1) =
\frac{(5a+2)(a+1)}{2}k^{a} + O(k^{a-1}).
$$
For big values of $k$ there exists a constant $C$ such that
$$
C (k^{1+a} (k+2a) - (k-1)^{1+a} (k+3a+1)) \ge a.
$$
But we can not choose a constant $C$ if $k^{1+a} (k+2a) \le (k-1)^{1+a} (k+3a+1)$.
There is a finite number of $k$ with $\frac{(5a+2)(a+1)}{2}k^{a} + O(k^{a-1}) \le 0$.
For such $k$ we want to prove that
$\MExpect N_n(1,k) = c(1,k) \,\left(n + O\left(f(k)\right)\right)$ with some function $f(k)$.
Using the method above we obtain the same inequalities:
$$
f(k) (k+2a) \ge
a
+  f(k-1) (k+3a+1)  ,
$$
$$
 f(k) (k+2a)\ge
a + C_1 k^{1+a} (k+3a+1).
$$
There exists a function $f$ such that the inequalities hold.
This completes the proof for $\MExpect N_n(1,k)$.

Consider the case $l>1$.
Assume that we already proved that
$$\MExpect N_n(i,j) = c(i,j) \, \left(n + \theta\left(C(i+j)^{1+a}\right)\right)$$
for all $i$ and $j$, such that $i<l$, $j \le k$ or $i\le l$ , $j < k$.

We use the following equality, which is similar to \eqref{M1}:
\begin{multline}\label{M2}
\MExpect N_{i+1} (l,k) =
\MExpect N_i(l,k) \left(1 - \frac {l(1+a)+k+a-1} {(a+1)i + a} \right) +\\
 + \frac {(l-2+a)\,\MExpect N_i(l-1,k)} {(a+1)i + a}
+ \frac {(k+al-1) \MExpect N_i(l,k-1)} {(a+1)i + a}.
\end{multline}

Note that if we have at least one vertex with first degree $l$ and second degree $k$ in $H_{a,1}^i$ (without a loop),
then we have at least $l+k-1$ edges in this graph.
Therefore $\MExpect N_{i}(l,k) = 0$ when $i < l+k-1$.
Consider the case $i = l+k-1$.
It is sufficient to prove that
$$
\MExpect N_{l+k-1} (l,k) \le C c(l,k) (l+k)
$$
with some $C$.
For any finite number of small $l$ and $k$ we can easily find a constant $C$ such that
$$
\MExpect N_{l+k-1} (l,k) \le C c(l,k) (l+k).
$$

Using \eqref{M2}, we get
$$
\MExpect N_{l+k-1} (l,k) =
\frac {(l-2+a)\,\MExpect N_{l+k-2}(l-1,k)} {(a+1)(l+k-2) + a}
+ \frac {(k+al-1) \MExpect N_{l+k-2}(l,k-1)} {(a+1)(l+k-2) + a} \le
$$
$$
\le Cc(l-1,k) \frac {(l-2+a)\left(l+k-1\right)} {(a+1)(l+k-2) + a}
+ Cc(l,k-1) \frac {(k+al-1)\left(l+k-1\right)} {(a+1)(l+k-2) + a} \le
$$
$$
\le Cc(l-1,k) \frac{(l-2+a)(l+k)}{l+al+k+2a}+
 C c(l,k-1) \frac{(al+k-1)(l+k)}{l+al+k+2a}.
$$
The last inequality holds if $k$ is big enough.

We also need to consider a finite number of small $k$.
First we show that for any finite number of small $k$ we have
$$
c(l,k) = \Omega \left( \frac{l^{k-4+\frac{a^2}{a+1}}}{(1+a)^l} \right).
$$
Indeed, from the recurrent relation we obtain
$$
c(l,1)=c(l-1,1)\frac{l-2+a}{(a+1)(l+1+a/(a+1))}.
$$
Therefore
$$
c(l,1)= \Omega\left( \frac{\Gamma(l-1+a)}{(a+1)^{l}\Gamma(l+2+a/(a+1))} \right) =
\Omega \left( \frac{l^{-3+\frac{a^2}{a+1}}}{(1+a)^l} \right).
$$
Here we used Statement \ref{stat1} from Subsection \ref{ss42}.
For $k\ge 2$ we have
$$
c(l,k)=c(l,k-1)\frac{al+k-1}{l(1+a)+k+2a}+c(l-1,k)\frac{l-2+a}{l(1+a)+k+2a}.
$$
It is sufficient to prove that there exists a positive function $f(k)$ such that
$$
f(k)(l(1+a)+k+2a)l^{k-4+\frac{a^2}{a+1}} \le
f(k-1)(al+k-1)l^{k-5+\frac{a^2}{a+1}} +
f(k)(l-2+a)(a+1)(l-1)^{k-4+\frac{a^2}{a+1}},
$$
$$
f(k)(l(1+a)+k+2a)\left(l^{k-4+\frac{a^2}{a+1}}-(l-1)^{k-4+\frac{a^2}{a+1}}\right)
+ f(k)(3a+k-a^2+2)(l-1)^{k-4+\frac{a^2}{a+1}}\le
$$
$$
\le f(k-1)(al+k-1)l^{k-5+\frac{a^2}{a+1}}.
$$
The last inequality holds for some positive function $f(k)$.

So we want to prove that
$$
\MExpect N_{l+k-1} (l,k) = O\left(  \frac{l^{k-3+\frac{a^2}{a+1}}}{(1+a)^l}  \right).
$$
Suppose that we have a graph on $l+k-1$ vertices and a vertex $t$ has first degree $l$ and second degree $k$.
Then one edge from this vertex goes to the vertex $1$, $l-1$ vertices send edges to $t$, and $k-2$ vertices send edges to the neighbors of $t$.
There are $\begin{pmatrix} l+k-2 \\ k-2 \end{pmatrix}$ ways to choose our vertex $t$ and its neighbors.
In each case the probability of these neighbors to send edges to the vertex $t$ equals
$$
O\left(\frac{(a(a+1)\dots(a+l-2))}{(3a+2)(4a+3)\dots(l(a+1)+a)}\right) =
O\left(\frac{\Gamma(a+l-1)}{(a+1)^l\Gamma(l+1+a/(a+1))}\right)
= O\left(\frac{l^{-2+\frac{a^2}{a+1}}}{(a+1)^l}\right),
$$
so
$$
\MExpect N_{l+k-1} (l,k) = O\left(\frac{l^{k-4+\frac{a^2}{a+1}}}{(a+1)^l}\right).
$$
This concludes the case $i = l+k-1$.

For $i \ge l + k-1$ we have
$$
\MExpect N_{i+1} (l,k) =
\MExpect N_i(l,k) \left(1 - \frac {l(1+a)+k+a-1} {(a+1)i + a} \right) + \frac {(l-2+a)\,\MExpect N_i(l-1,k)} {(a+1)i + a}
+
$$
$$
 +\frac {(k+al-1) \MExpect N_i(l,k-1)} {(a+1)i + a} = c(l,k) \left(i + \theta\left(C(l+k)^{1+a}\right)\right) \left(1 - \frac {l(1+a)+k+a-1} {(a+1)i + a}\right) +
$$
$$
+ c(l-1,k) \left(i + \theta\left(C(l+k-1)^{1+a}\right)\right) \frac {(l-2+a)} {(a+1)i + a}
+  c(l,k-1) \left(i + \theta\left(C(l+k-1)^{1+a}\right)\right) \frac {(k+al-1)} {(a+1)i + a} =
$$
$$
= c(l,k) i
- c(l,k) i \frac {l(1+a)+k+a-1} {(a+1)i + a}
+ c(l,k)\theta\left(C(l+k)^{1+a}\right) \left(1 - \frac {l(1+a)+k+a-1} {(a+1)i + a}\right) +
$$
$$
+ c(l-1,k) i \frac {(l-2+a)} {(a+1)i + a}
+  c(l,k-1) i \frac {(k+al-1)} {(a+1)i + a} +
$$
$$
+ c(l-1,k) \theta\left(C(l+k-1)^{1+a}\right)\frac {(l-2+a)} {(a+1)i + a}
+  c(l,k-1)  \theta\left(C(l+k-1)^{1+a}\right)\frac {(k+al-1)} {(a+1)i + a} =
$$
$$
= c(l,k) (i+1)
- c(l,k) \frac {il(1+a)+ik+2ia+ a} {(a+1)i + a}
+ c(l-1,k) i \frac {(l-2+a)} {(a+1)i + a}
+
$$
$$
+c(l,k-1) i \frac {(k+al-1)} {(a+1)i + a}+ c(l,k)\theta\left(C(l+k)^{1+a}\right) \left(1 - \frac {l(1+a)+k+a-1} {(a+1)i + a}\right) +
$$
$$
+ c(l-1,k) \theta\left(C(l+k-1)^{1+a}\right)\frac {(l-2+a)} {(a+1)i + a}
+  c(l,k-1)  \theta\left(C(l+k-1)^{1+a}\right)\frac {(k+al-1)} {(a+1)i + a} =
$$
$$
= c(l,k) (i+1)
-  \frac {a(k+al-1)c(l,k-1)} {((a+1)i + a)(l(1+a)+k+2a)}
-  \frac {a(l-2+a)c(l-1,k)} {((a+1)i + a)(l(1+a)+k+2a)} +
$$
$$
+ c(l,k)\theta\left(C(l+k)^{1+a}\right) \left(1 - \frac {l(1+a)+k+a-1} {(a+1)i + a}\right) +
$$
$$
+ c(l-1,k) \theta\left(C(l+k-1)^{1+a}\right)\frac {(l-2+a)} {(a+1)i + a}
+  c(l,k-1)  \theta\left(C(l+k-1)^{1+a}\right)\frac {(k+al-1)} {(a+1)i + a}.
$$
We want to prove the following inequality:
$$
C c(l,k) \left((l+k)^{1+a}\right) \frac {l(1+a)+k+a-1} {(a+1)i + a} \ge
$$
$$
\ge \frac {a(k+al-1)c(l,k-1)} {((a+1)i + a)(l(1+a)+k+2a)}
+ \frac {a(l-2+a)c(l-1,k)} {((a+1)i + a)(l(1+a)+k+2a)} +
$$
$$
+ C c(l-1,k) \left((l+k-1)^{1+a}\right)\frac {(l-2+a)} {(a+1)i + a}
+  C c(l,k-1)  \left((l+k-1)^{1+a}\right)\frac {(k+al-1)} {(a+1)i + a}.
$$
It is sufficient to show that the following inequalities hold
$$
C c(l,k-1) (l+k)^{1+a} \frac {(l(1+a)+k+a-1)(k+al-1)} {((a+1)i + a)(l(1+a)+k+2a)} \ge
$$
$$
\ge \frac {a(k+al-1)c(l,k-1)} {((a+1)i + a)(l(1+a)+k+2a)}
+  C c(l,k-1)  (l+k-1)^{1+a}\frac {(k+al-1)} {(a+1)i + a}
$$
and
$$
C c(l-1,k)(l+k)^{1+a}\frac {(l(1+a)+k+a-1)(l-2+a)} {((a+1)i + a)(l(1+a)+k+2a)} \ge
$$
$$
\ge  \frac {a(l-2+a)c(l-1,k)} {((a+1)i + a)(l(1+a)+k+2a)}
+ C c(l-1,k)(l+k-1)^{1+a}\frac {(l-2+a)} {(a+1)i + a}.
$$
In other words
$$
C(l+k)^{1+a}(l(1+a)+k+a-1)(k+al-1) \ge
$$
$$
\ge a(k+al-1)
+  C (l+k-1)^{1+a}(k+al-1)(l(1+a)+k+2a)
$$
and
$$
C (l+k)^{1+a} (l(1+a)+k+a-1)(l-2+a) \ge
$$
$$
\ge  a(l-2+a)
+ C  (l+k-1)^{1+a}(l-2+a)(l(1+a)+k+2a).
$$
To prove both inequalities we make the following transformations:
$$
(l+k)^{1+a}(l(1+a)+k+a-1) -
(l+k-1)^{1+a}(l(1+a)+k+2a) =
$$
$$
= (l+k)^{1+a}(l(1+a)+k+a-1)
- ((l+k)^{1+a} - (1+a)(l+k)^{a} + \frac{a(1+a)}{2}(l+k)^{a-1} +
$$
$$
+ O\left((l+k)^{a-2})\right)(l(1+a)+k+2a) = -(l+k)^{1+a}(1+a) + (1+a)(l+k)^{a}(l(1+a)+k+2a) -
$$
$$
- \frac{a(1+a)}{2}(l+k)^{a-1}(l(1+a)+k+2a)
+ O\left((l+k)^{a-2}\right)(l(1+a)+k+2a) =
$$
$$
= (l+k)^{a-1}(1+a)
\left( al^2 + alk + 2al + 2ak
- \frac{a(1+a)}{2}l - \frac{a}{2}k - \frac{2a^2}{2} \right)
+ O\left((l+k)^{a-2}\right)(l(1+a)+k+2a) =
$$
$$
= (l+k)^{a-1}(1+a)
\left( al^2 + alk + \frac{3}{2}al + \frac{3}{2}ak
- \frac{1}{2}a^2l - a^2\right)
+ O\left((l+k)^{a-2}\right)(l(1+a)+k+2a).
$$
If $l$ or $k$ is large enough, then there exists a constant $C$ such that
$$
C(l+k)^{a-1}(1+a)
\left( al^2 + alk + \frac{3}{2}al + \frac{3}{2}ak
- \frac{1}{2}a^2l - a^2\right)
+ O\left((l+k)^{a-2}\right)(l(1+a)+k+2a) \ge  a.
$$

Finally, we need to consider the finite number of small $l$ and $k$.
We want to find some function $f(l,k)$ such that
$$
f(l,k) c(l,k) \left((l+k)^{1+a}\right) \frac {l(1+a)+k+a-1} {(a+1)i + a} \ge
$$
$$
\ge \frac {a(k+al-1)c(l,k-1)} {((a+1)i + a)(l(1+a)+k+2a)}
+ \frac {a(l-2+a)c(l-1,k)} {((a+1)i + a)(l(1+a)+k+2a)} +
$$
$$
+ f(l-1,k) c(l-1,k) \frac {(l-2+a)} {(a+1)i + a}
+  f(l,k-1) c(l,k-1) \frac {(k+al-1)} {(a+1)i + a}.
$$
Such function $f(l,k)$ exists.
This concludes the proof of Theorem \ref{th4}.

\subsection{Proof of Theorem \ref{exp}}\label{ss42}

In this proof we shall use the following statement.
\begin{Stat}\label{stat1}
For $t>0$ and fixed $a>0$
$$
\frac{\Gamma(t+a)}{\Gamma(t)} = t^a\left(1 + O(1/t)\right).
$$
\end{Stat}

\begin{Proof}
From Stirling's formula we obtain
$$
\frac{\Gamma(t+a)}{\Gamma(t)} =
\sqrt{\frac{t}{t+a}}\frac{(t+a)^{a}}{e^{a}}\left(\frac{t+a}{t}\right)^t\left(1 + 1/t\right).
$$
It is easy to check that
$$
t \ln \left(1 + \frac{a}{t}\right) = a + O(1/t).
$$
So
$$
\left(1 + \frac{a}{t}\right)^t = e^a (1 + O(1/t)).
$$
We obtain
$$
\frac{\Gamma(t+a)}{\Gamma(t)}
= \sqrt{\frac{t}{t+a}}(t+a)^{a}\left(1 + O(1/t)\right)
= t^{a}\left(1 + O(1/t)\right).
$$
\end{Proof}

\subsubsection{Estimation of $c(1,k)$}

\begin{Lem}
$$
c(1,k) = \frac{\Gamma(2a+1)(1 + O(1/k))}{\Gamma(a)k^{a+1}}.
$$
\end{Lem}

\begin{Proof}
As we know
$$
c(k) = \frac {\B(k-1+a,a+2)} {\B(a,a+1)} =
\frac{(a+1)\Gamma(2a+1)\Gamma(k-1+a)}{\Gamma(a)\Gamma(k+1+2a)}.
$$
Using the recurrent relation
$$
c(1,k) = c(1,k-1)\frac{a+k-1}{k+3a+1}+c(k)\frac{a+k-1}{k+3a+1}
$$
we obtain
$$
c(1,k) = \sum_{j=1}^{k} \frac{c(j)(a+j-1)\ldots(a+k-1)}{(j+3a+1)\ldots(k+3a+1)} =
$$
$$
= \frac{(a+1)\Gamma(2a+1)}{\Gamma(a)}\sum_{j=1}^{k} \frac{\Gamma(j-1+a)(a+j-1)\ldots(a+k-1)}{\Gamma(j+1+2a)(j+3a+1)\ldots(k+3a+1)} =
$$
$$
= \frac{(a+1)\Gamma(2a+1)\Gamma(a+k)}{\Gamma(a)\Gamma(k+3a+2)}\sum_{j=1}^{k} \frac{\Gamma(j+3a+1)}{\Gamma(j+1+2a)} =
$$
$$
= \frac{(a+1)\Gamma(2a+1)\Gamma(a+k)}{\Gamma(a)\Gamma(k+3a+2)}\sum_{j=1}^{k} j^a (1 + O(1/j))
= \frac{\Gamma(2a+1)k^{a+1}(1 + O(1/k))}{\Gamma(a)k^{2a+2}} =
$$
$$
= \frac{\Gamma(2a+1)(1 + O(1/k))}{\Gamma(a)k^{a+1}}.
$$
\end{Proof}

\subsubsection{Sum of $c(l,k)$}\label{ss44}

We want to estimate the sum
$\sum_{l=1}^{\infty} c(l,k)$.
First let us prove that the series
$\sum_{l=1}^{\infty}l^{N}c(l,k)$ converges for all $N$ and $k$.

The inequality
$$
c(l,k) \le \tilde C \frac{p^{k}}{(1+q)^{l}}
$$
holds for any $p > 1$ and $q = \min\{a,1\}\frac{(p-1)}{p}$.
Here we choose $\tilde C$ so that $\tilde C \frac{p^{k}}{1+ap} \ge c(1,k)$ for any $k$.
We need to prove that
$$
\frac{p^{k}}{(1+q)^{l}} (l+al+k+2a) \ge
\frac{p^{k-1}}{(1+q)^{l}}(al+k-1) + \frac{p^{k}}{(1+q)^{l-1}}(l-2+a),
$$
We make some transformations:
$$
p (l+al+k+2a) \ge
(al+k-1) + p(1+q)(l-2+a),
$$
$$
p (al+k+a+2)  \ge
(al+k-1) + pq(l-2+a),
$$
$$
al+k+a+2  \ge
\min\{a,1\}(l-2+a).
$$
The last inequality holds.
Therefore $c(l,k) \le \tilde C \frac{p^{k}}{(1+q)^{l}}$ and
$\sum_{l=1}^{\infty}l^{N}c(l,k)$ converges.

For $l \ge 2$ and any $c \ge 0$ we have
$$
c(l,k)(l(1+a)+k+2a) \frac{\Gamma(l+a+c)}{\Gamma(l+a-1)}
= c(l,k-1) \frac{\Gamma(l+a+c)}{\Gamma(l+a-1)} (al+k-1) +
c(l-1,k)(l-2+a) \frac{\Gamma(l+a+c)}{\Gamma(l+a-1)} .
$$
Therefore,
$$
\sum_{l = 2}^{\infty} c(l,k)(l(1+a)+k+2a) \frac{\Gamma(l+a+c)}{\Gamma(l+a-1)}
= \sum_{l = 2}^{\infty} c(l,k-1)(al+k-1) \frac{\Gamma(l+a+c)}{\Gamma(l+a-1)}
+ \sum_{l = 1}^{\infty} c(l,k)\frac{\Gamma(l+a+c+1)}{\Gamma(l+a-1)} ,
$$
$$
\sum_{l = 2}^{\infty} c(l,k)(al+k+a-c) \frac{\Gamma(l+a+c)}{\Gamma(l+a-1)}
= \sum_{l = 2}^{\infty} c(l,k-1)(al+k-1) \frac{\Gamma(l+a+c)}{\Gamma(l+a-1)}
+ c(1,k)\frac{\Gamma(a+c+2)}{\Gamma(a)}.
$$
Consider the function
$$
f_c(k) = \frac{\Gamma(k+a-c)}{\Gamma(k)} = k^{a-c}(1+O(1/k)).
$$
It is easy to see that
$$
\frac{f_c(k+1)}{f_c(k)} = 1 + \frac{a-c}{k}.
$$
We have
$$
\sum_{j=1}^k \sum_{l = 2}^{\infty} c(l,j)(al+j+a-c) \frac{\Gamma(l+a+c)}{\Gamma(l+a-1)} f_c(j) =
$$
$$
= \sum_{j=1}^k \sum_{l = 2}^{\infty} c(l,j-1)(al+j-1) \frac{\Gamma(l+a+c)}{\Gamma(l+a-1)} f_c(j)
+ \sum_{j=1}^k c(1,j)\frac{\Gamma(a+c+2)}{\Gamma(a)} f_c(j) ,
$$
$$
\sum_{j=1}^k \sum_{l = 2}^{\infty} c(l,j)(al+j+a-c) \frac{\Gamma(l+a+c)}{\Gamma(l+a-1)} f_c(j) =
$$
$$
= \sum_{j=1}^{k-1} \sum_{l = 2}^{\infty} c(l,j)(al+j) \frac{\Gamma(l+a+c)}{\Gamma(l+a-1)} f_c(j)\left(1 + \frac{a-c}{j} \right)
+ \sum_{j=1}^k  c(1,j)\frac{\Gamma(a+c+2)}{\Gamma(a)} f_c(j) ,
$$
$$
\sum_{l = 2}^{\infty} c(l,k)(al+k+a-c) \frac{\Gamma(l+a+c)}{\Gamma(l+a-1)} f_c(k) =
$$
$$
= \sum_{j=1}^{k-1} \sum_{l = 2}^{\infty} c(l,j) \frac{\Gamma(l+a+c)}{\Gamma(l+a-1)}\frac{al(a-c)}{j} f_c(j)
+ \sum_{j=1}^k  c(1,j)\frac{\Gamma(a+c+2)}{\Gamma(a)} f_c(j).
$$
If $c \ge a$ then taking into account Lemma 8 and the above-mentioned asymptotics for $f_c(j) $ we have
$$
\sum_{l = 2}^{\infty} c(l,k)(al+k+a-c) \frac{\Gamma(l+a+c)}{\Gamma(l+a-1)} f_c(k) \le
\sum_{j=1}^k  c(1,j)\frac{\Gamma(a+c+2)}{\Gamma(a)} f_c(j) = O(1).
$$
Hence,
$$
\sum_{l = 2}^{\infty} c(l,k) \frac{\Gamma(l+a+c)}{\Gamma(l+a-1)} f_c(k) = O(1/k).
$$
We want to prove that for any $0 \le c < a+1$ the following equality holds:
\begin{equation}\label{sum}
\sum_{l = 2}^{\infty} c(l,k) \frac{\Gamma(l+a+c)}{\Gamma(l+a-1)} f_c(k) = O\left(\frac{(\ln k)^{\lceil a-c\rceil }}{k}\right).
\end{equation}
We have already proved this statement for $a \le c < a+1$.

Suppose that for $c' \ge 1$ we have
\begin{equation*}\label{eq11}
\sum_{l = 2}^{\infty} c(l,k) \frac{\Gamma(l+a+c')}{\Gamma(l+a-1)} f_{c'}(k)
= O\left(\frac{(\ln k)^{\lceil a-c'\rceil }}{k}\right).
\end{equation*}
Then
$$
\sum_{l = 2}^{\infty} c(l,k)(al+k+a-c'+1) \frac{\Gamma(l+a+c'-1)}{\Gamma(l+a-1)} f_{c'-1}(k) =
$$
$$
= \sum_{j=1}^{k-1} \sum_{l = 2}^{\infty} c(l,j) \frac{\Gamma(l+a+c'-1)}{\Gamma(l+a-1)}\frac{al(a-c'+1)}{j}(j+a-c')f_{c'}(j)
+ \sum_{j=1}^k  c(1,j)\frac{\Gamma(a+c'+1)}{\Gamma(a)} f_{c'-1}(j) =
$$
$$
= O\left( \sum_{j=1}^{k-1} \frac{(\ln k)^{\lceil a-c'\rceil }}{j} \right)
= O\left((\ln k)^{\lceil a-c'+1\rceil }\right).
$$
We proved \eqref{sum}. In particular,
\begin{equation}\label{sum'}
\sum_{l = 2}^{\infty} c(l,k) \frac{\Gamma(l+a)}{\Gamma(l+a-1)} f_0(k)
= \sum_{l = 2}^{\infty} c(l,k) (l+a-1) f_0(k)
= O\left(\frac{(\ln k)^{\lceil a\rceil }}{k}\right).
\end{equation}

Put $x_k = \sum_{l = 2}^{\infty} c(l,k)$.
For $l \ge 2$
$$
c(l,k)(l(1+a)+k+2a) = c(l,k-1)(al+k-1)+c(l-1,k)(l-2+a).
$$
So
$$
\sum_{l = 2}^{\infty} c(l,k)(l(1+a)+k+2a)
= \sum_{l = 2}^{\infty} c(l,k-1)(al+k-1) + \sum_{l = 1}^{\infty} c(l,k)(l-1+a),
$$
$$
\sum_{l = 2}^{\infty} c(l,k)(al+k+a+1) = \sum_{l = 2}^{\infty} c(l,k-1)(al+k-1) + a c(1,k),
$$
$$
(k+a+1)x_k = (k-1)x_{k-1} + a c(1,k)
+ a \sum_{l = 2}^{\infty}l(c(l,k-1) - c(l,k)).
$$
We have
$$
(k+a+1)f_{-1}(k)x_k = (k-1)f_{-1}(k)x_{k-1} + af_{-1}(k)c(1,k)
+ a f_{-1}(k)\sum_{l = 2}^{\infty}l(c(l,k-1) - c(l,k)),
$$
$$
\sum_{j=1}^k (j+a+1)f_{-1}(j)x_j =
\sum_{j=1}^{k-1}f_{-1}(j)(j + a + 1) x_{j}
+ \sum_{j=1}^k a f_{-1}(j)c(1,j)
+ \sum_{j=1}^k a f_{-1}(j)\sum_{l = 2}^{\infty}l(c(l,j-1) - c(l,j)),
$$
$$
(k+a+1)f_{-1}(k)x_k =
a \sum_{j=1}^k f_{-1}(j)c(1,j)
+ a \sum_{j=1}^k f_{-1}(j)\sum_{l = 2}^{\infty}l(c(l,j-1) - c(l,j)).
$$
$$
f_{-1}(k)(k+a+1)x_k =
a \sum_{j=1}^k f_{-1}(j)c(1,j) +
a(a+1) \sum_{j=1}^{k-1}\frac{ f_{-1}(j)}{j} \sum_{l = 2}^{\infty}l c(l,j)
- a f_{-1}(k) \sum_{l = 2}^{\infty}l c(l,k) =
$$
$$
= a \sum_{j=1}^k j^{a+1} \frac{\Gamma(2a+1)}{\Gamma(a)j^{a+1}}(1 + O(1/j)) +
\sum_{j=1}^{k-1}    O\left(\frac{(\ln j)^{\lceil a\rceil }}{j}\right)
+ O\left((\ln k)^{\lceil a\rceil} \right) =
$$
$$
= a k \frac{\Gamma(2a+1)}{\Gamma(a)}  +
\sum_{j=1}^{k-1}   O\left(\frac{(\ln j)^{\lceil a\rceil }}{j}\right)
+ O\left((\ln k)^{\lceil a\rceil} \right) =
a k \frac{\Gamma(2a+1)}{\Gamma(a)} \left(1 +  O\left(\frac{(\ln k)^{\lceil a+1\rceil }}{k}\right) \right).
$$
Here we used \eqref{sum'} and Lemma 8. We obtain
$$
x_k = \frac{a\Gamma(2a+1)}{\Gamma(a)k^{a+1}} \left(1 +  O\left(\frac{(\ln k)^{\lceil a+1\rceil }}{k}\right) \right)
$$
and
$$
\sum_{l=1}^{\infty} c(l,k) = c(1,k) + x_k =
\frac{(a+1)\Gamma(2a+1)}{\Gamma(a)k^{a+1}} \left(1 +  O\left(\frac{(\ln k)^{\lceil a+1\rceil }}{k}\right) \right).
$$

\subsubsection{Estimation of $\MExpect Y_n(k)$ }

Note that
$$
\sum_{l\ge1}\sum_{j\ge k} \MExpect N_{i+1} (l,j) =
\sum_{l\ge1}\sum_{j\ge k}\MExpect N_i(l,j)
+ \sum_{l\ge1}\frac {(al+k-1) \MExpect N_i(l,k-1)} {(a + 1)i + a}
+ \sum_{j\ge k}\frac {(j-1+a)\,M_i^1(j)} {(a + 1)i + a}.
$$
Therefore we obtain
$$
\sum_{l\ge1}\sum_{j\ge k} \MExpect N_{n} (l,j) =
\sum_{i=1}^{n-1} \sum_{l\ge1}\frac {(al+k-1) \MExpect N_i(l,k-1)} {(a + 1)i + a}
+ \sum_{i=1}^{n-1} \sum_{j\ge k}\frac {(j-1+a)\,M_i^1(j)} {(a + 1)i + a}.
$$

Let us estimate the sum
$$
\sum_{i=1}^{n-1} \sum_{j\ge k}\frac {(j-1+a)\,M_i^1(j)} {(a + 1)i + a}.
$$
First we compute
$$
F_t(k) = \sum_{j\ge k}(j-1+a)\,M_t^1(j).
$$
Let us prove by induction on $k$ that
$$
F_n(k)
=\frac{(a+1)\Gamma(2a+1)\Gamma(k+a)}{\Gamma(a+1)\Gamma(k+2a)}n \left(1+\theta\left(\frac{C(k-1)^{1+a}}{n}\right)\right)
$$
with some constant $C$.
For $k = 1$ and $k=2$ we have
$$
F_n(1) = \sum_{j\ge 1}(j-1+a)\,M_n^1(j) = n(1+a),
$$
$$
F_n(2) = F_n(1) - a\,M_n^1(1) = n(1+a) - a n \frac{(1+a)}{(2a+1)} +  O\left(1\right)
= \frac{(1+a)^2 n}{2a+1}\left(1 + O\left(1/n\right)\right).
$$
For $k \ge 3$ we have
$$
M_{i+1}^1(j) =
M_{i}^1(j)\left(1 - \frac{j-1+a}{(a + 1)i + a}\right) +
M_{i}^1(j-1)\frac{j-2+a}{(a + 1)i + a}.
$$
We multiply this equality by $(j-1+a)$ and sum over all $j \ge k$:
$$
F_{i+1}(k) = \sum_{j\ge k} (j-1+a) \, M_{i+1}^1(j) =
$$
$$
= \sum_{j\ge k} (j-1+a) M_{i}^1(j)
 -\sum_{j\ge k} M_{i}^1(j) \frac{(j-1+a)(j-1+a)}{(a + 1)i + a} +
\sum_{j\ge k-1} M_{i}^1(j)\frac{(j+a)(j-1+a)}{(a + 1)i + a} =
$$
$$
= F_i(k)
 + \sum_{j\ge k} M_{i}^1(j) \frac{(j-1+a)}{(a + 1)i + a} +
M_{i}^1(k-1)\frac{(k-1+a)(k-2+a)}{(a + 1)i + a} =
$$
$$
= F_i(k)\left(1 + \frac{1}{(a + 1)i + a}\right) +
\left(F_i(k-1) - F_i(k)\right)\frac{(k-1+a)}{(a + 1)i + a} =
$$
$$
= F_i(k)\left(1 - \frac{k-2+a}{(a + 1)i + a}\right) +
F_i(k-1)\frac{(k-1+a)}{(a + 1)i + a}.
$$

Note that for $i+1 < k-1$ we have $F_{i+1}(k) = 0$.
Consider $i+1 \ge k-1$.
Using the inductive assumption we get
$$
F_{i+1}(k) =
\frac{(a+1)\Gamma(2a+1)\Gamma(k+a)}{\Gamma(a+1)\Gamma(k+2a)}i
\left(1 - \frac{k-2+a}{(a + 1)i + a}\right)\left(1+\theta\left(\frac{C(k-1)^{1+a}}{i}\right)\right) +
$$
$$
+ \frac{(a+1)\Gamma(2a+1)\Gamma(k-1+a)i}{\Gamma(a+1)\Gamma(k-1+2a)}\frac{(k-1+a)}{(a + 1)i + a}
\left(1+\theta\left(\frac{C(k-2)^{1+a}}{i}\right)\right) =
$$
$$
= \frac{(a+1)\Gamma(2a+1)\Gamma(k+a)}{\Gamma(a+1)\Gamma(k+2a)}
\Biggl( i+1 -\frac{a}{(a + 1)i + a}
+ i \left(1 - \frac{k-2+a}{(a + 1)i + a}\right)\theta\left(\frac{C(k-1)^{1+a}}{i}\right) +
$$
$$
+ \frac{(k-1+2a)i}{(a + 1)i + a}
\theta\left(\frac{C(k-2)^{1+a}}{i}\right)
 \Biggr).
$$
And we need to show that for some constant $C$
$$
\frac{a}{(a + 1)i + a}
+ \frac{(k-1+2a)}{(a + 1)i + a}C(k-2)^{1+a}
\le \frac{k-2+a}{(a + 1)i + a} C(k-1)^{1+a}.
$$
This inequality holds for sufficiently large $C$.

We have
$$
\sum_{i=1}^{n-1} \frac{F_i(k)}{(a+1)i+a} =
\sum_{i=1}^{n-1} \frac{\Gamma(2a+1)\Gamma(k+a)}{\Gamma(a+1)\Gamma(k+2a)} \left(1+O\left(\frac{k^{1+a}}{i}\right)\right) =
$$
$$
= \frac{\Gamma(2a+1)\Gamma(k+a)}{\Gamma(a+1)\Gamma(k+2a)}n \left(1+O\left(\frac{k^{1+a}}{n}\right)\right).
$$

Let us estimate the sum
$$
\sum_{i=1}^{n-1} \sum_{l\ge1}\frac {(al+k-1) \MExpect N_i(l,k-1)} {(a + 1)i + a}.
$$
We start with the sum
$$
\sum_{l\ge1} (al+k-1) \MExpect N_i(l,k-1).
$$
It is easy to see that
$$
\MExpect N_i(l,k) = O(c(l,k)i).
$$
To verify this, one can follow the proof of Theorem \ref{th4} and make sure that it works for the inequality
$$\MExpect N_i(l,k) < \tilde{C} c(l,k)((a+1)i+a)$$ with some constant $\tilde{C}$ -- note that the analog of Lemma \ref{Gr} is also needed.

Therefore
$$
\sum_{l\ge1} (al-1) \MExpect N_i(l,k-1) =
O\left(\sum_{l\ge1}(al-1)c(l,k-1)i\right) =
O\left(\frac{(\ln k)^{\lceil a\rceil }i}{k^{a+1}}\right).
$$
Using \eqref{sum}
we obtain
$$
\sum_{l\ge1}k \MExpect N_i(l,k-1) =
\sum_{l\ge1}k c(l,k-1)i \left(1+ O\left((l+k)^{1+a}/i\right)\right) =
$$
$$
= \frac{(a+1)\Gamma(2a+1)}{\Gamma(a)k^{a}} i\left(1 +  O\left(\frac{(\ln k)^{\lceil a+1\rceil }}{k}\right)
+ O\left(\frac{k^{1+a}}{i}\right) \right).
$$
Here we used the following estimate:
$$
\sum_{l\ge1}k c(l,k-1)(l+k)^{1+a} =
O\left(\sum_{l=1}^{k}k^{2+a}c(l,k-1) + \sum_{l \ge k}k c(l,k-1)l^{1+a}\right) =
O(k)+O(1) = O(k).
$$
So we have
$$
\sum_{i=1}^{n-1} \sum_{l\ge1}\frac {(al+k-1) \MExpect N_i(l,k-1)} {(a + 1)i + a} =
 \frac{\Gamma(2a+1)}{\Gamma(a)k^{a}} n \left(1 +  O\left(\frac{(\ln k)^{\lceil a+1\rceil }}{k}\right)
 + O\left(\frac{k^{1+a}}{n}\right) \right).
$$
Hence
$$
\sum_{l\ge1}\sum_{j\ge k} \MExpect N_{n} (l,j) =
\frac{a\Gamma(2a+1)}{\Gamma(a+1)k^{a}} n \left(1 +  O\left(\frac{(\ln k)^{\lceil a+1\rceil }}{k}\right) \right)
+ \frac{\Gamma(2a+1)\Gamma(k+a)}{\Gamma(a+1)\Gamma(k+2a)}n \left(1+O\left(\frac{k^{1+a}}{n}\right)\right) =
$$
$$
= \frac{(a+1)\Gamma(2a+1)}{\Gamma(a+1)k^{a}} n \left(1 +  O\left(\frac{(\ln k)^{\lceil a+1\rceil }}{k}\right)
+O\left(\frac{k^{1+a}}{n}\right) \right).
$$

Consider vertices with loops.
For $k=0$, we have
$$
\sum_{l\ge1}\sum_{j\ge 0} \MExpect P_{n} (l,j) =
\sum_{i=1}^n \frac{a}{(1+a)i-1} = O(\ln n).
$$
For $k\ge2$, we have
$$
\sum_{l\ge1}\sum_{j\ge k} \MExpect P_{i+1} (l,j) =
\sum_{l\ge1}\sum_{j\ge k}\MExpect P_i(l,j)
+ \sum_{l\ge1}\frac {(al+k-2a-1) \MExpect P_i(l,k-1)} {(a + 1)i + a}.
$$
Therefore, we obtain
$$
\sum_{l\ge1}\sum_{j\ge k} \MExpect P_{n} (l,j) =
\sum_{i=1}^{n-1} \sum_{l\ge1}  \frac {(al+k-2a-1) \MExpect P_i(l,k-1)} {(a + 1)i + a} \le
\sum_{i=1}^{n-1} \sum_{l\ge1}  \frac {(al+k-2a-1) p(l,k-1)}{(a + 1)i}.
$$

From the recurrent relation for $p(l,k)$ it follows that
$$
p(l,k) = O\left(\frac{1}{l^2}\right)\,,
$$
and
$$
p(l,k) = O\left(\frac{k}{l^3}\right)\,.
$$
To obtain the second estimate consider $q(l,k) = p(l,k)/k$.
For $k\ge 1$ we have
$$
q(l,k)(l+al+k-1-a)=q(l,k-1)\frac{(k-1)(al+k-2a-1)}{k}+q(l-1,k)(l-2+a),
$$
$$
q(l,k)(l+al+k-1-a)-q(l,k-1)\left(al+k-2a-2-\frac{al-2a-1}{k}\right) = q(l-1,k)(l-2+a).
$$
Thus, $q(l,k) = O(q(l))$, where
$$
q(l)(l+a+1+(al-2a-1)) = q(l-1)(l-2+a).
$$
From this equality it follows that $q(l)=O\left(\frac{1}{l^3}\right)$.

We can estimate the following sum:
$$
\sum_{l\ge1}  \frac {(al+k-2a-1) p(l,k-1)}{(a + 1)} =
O\left(k\right).
$$
Hence
$$
\sum_{l\ge1}\sum_{j\ge k} \MExpect P_{n} (l,j) = O(k \ln n).
$$

Now we are ready to estimate $\MExpect Y_n(k)$:
$$
\MExpect Y_n(k) = \sum_{l\ge1}\sum_{j\ge k} \MExpect N_{n} (l,j) + \sum_{l\ge1}\sum_{j\ge k} \MExpect P_{n} (l,k) =
$$
$$
= \frac{(a+1)\Gamma(2a+1)}{\Gamma(a+1)k^{a}} n \left(1 +  O\left(\frac{(\ln k)^{\lceil a+1\rceil }}{k}\right)
+O\left(\frac{k^{1+a}}{n}\right) \right) + O(k \ln n) =
$$
$$
= \frac{(a+1)\Gamma(2a+1)}{\Gamma(a+1)k^{a}} n \left(1 +  O\left(\frac{(\ln k)^{\lceil a+1\rceil }}{k}\right)
+O\left(\frac{k^{1+a}}{n}\right) \right).
$$
This concludes the proof of Theorem \ref{exp}.
\subsection{Proof of Lemma \ref{lem2}}\label{ss43}

It is easy to see that $\MExpect P_n(0,k)= \MExpect P_n(1,k)=0$.
For all $k>0$ we have $\MExpect P_n(2,k)=0$.
For $k=0$ we have
$$
\MExpect P_n(2,0)
= \sum_{i=1}^n {\frac {a} {(a+1)i-1} \prod_{j=i+1}^{n} \frac {(1+a)j-2-a} {(1+a)j-1}}
= \sum_{i=1}^n \frac {a} {(a+1)i-1}
\frac {\Gamma\left(n-\frac{1}{a+1}\right)\Gamma\left(i+\frac{a}{a+1}\right)} {\Gamma\left(n+\frac{a}{a+1}\right)\Gamma\left(i-\frac{1}{a+1}\right)} =
$$
$$
= \frac{1}{n}(1+ O(1/n)) \sum_{i=1}^n {\frac {ai} {(a+1)i-1} (1 + O(1/i)) } = O(1).
$$

The rest of the proof is by induction.
Consider $l\ge3$, $k\ge1$.
Assume that we already proved that
$\MExpect P_n(i,j) \le p(i,j)$
for all $i$ and $j$, such that $i<l$, $j \le k$ or $i\le l$ , $j < k$.
We use the following equality
\begin{multline}\label{P}
\MExpect P_{i+1} (l,k) =
\MExpect P_i(l,k)\left(1-\frac{l(a+1)+k-a-1}{(a+1)i+a}\right) + \MExpect P_i(l,k-1)\frac{al+k-2a-1}{(a+1)i+a}+\\
+\MExpect P_i(l-1,k)\frac{l-2+a}{(a+1)i+a}.
\end{multline}

Note that if we have at least one vertex with a loop, with first degree $l$ and second degree $k$ in the graph $H_{a,1}^i$,
then we have at least $l+k-1$ edges in this graph.
Therefore $\MExpect P_{i}(l,k) = 0$ if $i < l+k-1$.
Consider the case $i = l+k-1$.
Using \eqref{P}, we get (for $k \ge 1$)
$$
\MExpect P_{l+k-1} (l,k) =
\frac {(l-2+a)\,\MExpect P_{l+k-2}(l-1,k)} {(a+1)(l+k-2) + a}
+ \frac {(al+k-2a-1) \MExpect P_{l+k-2}(l,k-1)} {(a+1)(l+k-2) + a} \le
$$
$$
\le \frac {(l-2+a) p(l-1,k)} {(a+1)(l+k-2) + a}
+ \frac {(al+k-2a-1) p(l,k-1)} {(a+1)(l+k-2) + a} =
$$
$$
= \frac {(l+al+k-1-a) p(l,k)} {(a+1)(l+k-2) + a}
\le p(l,k).
$$
The last inequality holds for $k \ge 1/a$.
Consider the case $k < 1/a$.
As in proof of Theorem \ref{th4}, at first we estimate $p(l,k)$:
$$
p(l,k) = \Omega\left(\frac{l^{k+\frac{a^2}{a+1}}}{(1+a)^l}\right).
$$
For $k=0$ we have
$$
p(l,0) = p(l-1,0)\frac{l-2+a}{(1+a)(l-1-\frac{1}{a+1})}.
$$
Therefore,
$$
p(l,0) = \Omega\left(\frac{l^{\frac{a^2}{a+1}}}{(1+a)^l}\right).
$$
For $k\ge 1$ we have
$$
p(l,k)=p(l,k-1)\frac{al+k-2a-1}{l(1+a)+k-1-a}+p(l-1,k)\frac{l-2+a}{l(1+a)+k-1-a}.
$$
Again, it is sufficient to prove that there exists a positive function $f(k)$ such that for big $l$
$$
f(k)(l(1+a)+k-1-a)l^{k+\frac{a^2}{a+1}} \le
f(k-1)(al+k-2a-1)l^{k+\frac{a^2}{a+1}-1} +
f(k)(l-2+a)(a+1)(l-1)^{k+\frac{a^2}{a+1}},
$$
$$
f(k)(l(1+a)+k-1-a)\left(l^{k+\frac{a^2}{a+1}}-(l-1)^{k+\frac{a^2}{a+1}}\right)
+f(k)(k-a^2+1)(l-1)^{k+\frac{a^2}{a+1}} \le
$$
$$
\le f(k-1)(al+k-2a-1)l^{k+\frac{a^2}{a+1}-1}.
$$
The last inequality holds for some function $f(k)$.

We want to prove that
$$
\MExpect P_{l+k-1} (l,k) = O\left(\frac{l^{k+\frac{a^2}{a+1}}}{(1+a)^l}\right).
$$
There are $l^k$ possible graphs on $l+k-1$ vertices with some vertex of first degree $l$, second degree $k$, and without a loop. And this vertex is exactly the vertex $1$.
The probability of this vertex to be a vertex with first degree $l$ and second degree $k$ equals
$$
O\left(\frac{l^k((a+1)\dots(a+l-2))}{(a+2)\dots((l-1)(a+1)-1)}\right)
= O\left(\frac{l^{k+\frac{a^2}{a+1}}}{(a+1)^l}\right).
$$
This concludes the case $i = l+k-1$.

If $i \ge l + k - 1$, then
$$
\MExpect P_{i+1} (l,k) =
\MExpect P_i(l,k)\left(1-\frac{l(a+1)+k-a-1}{(a+1)i+a}\right)+
$$
$$
 + \MExpect P_i(l,k-1)\frac{al+k-2a-1}{(a+1)i+a}
+\MExpect P_i(l-1,k)\frac{l-2+a}{(a+1)i+a}.
$$
Using the recurrent relation for $p(l,k)$ and induction on $i$ it is easy to prove that
$\MExpect P_n(l,k) \le p(l,k)$.
This concludes the proof of Lemma \ref{lem2}.

\subsection{Proof of Theorem \ref{X_n}}\label{Corol}

We estimate the expectation of $X_n(k)$ as follows:
$$
\MExpect X_n(k) = \sum_{l=1}^{\infty} \MExpect N_n(l,k) + \sum_{l=1}^{\infty} \MExpect P_n(l,k) =
\sum_{l=1}^{\infty} c(l,k) n + O\left( \sum_{l=1}^{\infty} c(l,k) (l+k)^{1+a} \right) +
O\left( \sum_{l=1}^{\infty} p(l,k) \right) =
$$
$$
= \frac{(a+1)\Gamma(2a+1)n}{\Gamma(a)k^{a+1}} \left(1 +  O\left(\frac{(\ln k)^{\lceil a+1\rceil }}{k}\right) \right) +
O(1) + O(1) =
$$
$$
= \frac{(a+1)\Gamma(2a+1)n}{\Gamma(a)k^{a+1}} \left(1 +  O\left(\frac{(\ln k)^{\lceil a+1\rceil }}{k}\right) +
O\left(\frac{k^{1+a}}{n}\right) \right).
$$

\renewcommand{\refname}{References}

\end{document}